\documentclass[12pt]{article}
\usepackage[psamsfonts]{amsfonts}
\usepackage{amsmath,amssymb,amsthm}
\usepackage{graphicx}
\usepackage{color}
\usepackage{hyperref}
\usepackage{eepic}
\newtheorem{THM}{Theorem}[section]

\newtheorem{LEM}[THM]{Lemma}
\newtheorem{PRP}[THM]{Proposition}

\newcommand{\sN}{{\sss N}}

\newcommand{\code}[1]{}

\def\UseSection{%%
        \numberwithin{equation}{section}
	\theoremstyle{plain}% default theorem style 
        \newtheorem{theorem}    {Theorem}[section]
        \DefineTheorems % Use this to define other environments to be 
}

\def\DefineTheorems{%%
	
	\newtheorem{lemma}      [theorem] {Lemma}
	
	\newtheorem{prop}       [theorem] {Proposition}
	
	\newtheorem{cor}        [theorem] {Corollary}

	\theoremstyle{definition}% ``defn'' theorem style 
	\newtheorem{defn}       [theorem] {Definition}

	\theoremstyle{definition}% ``remark'' theorem style 

}

\newcommand{\bt}   {\begin{theorem}}
\newcommand{\et}   {\end  {theorem}}
\newcommand{\bl}   {\begin{lemma}}
\newcommand{\el}   {\end  {lemma}}
\newcommand{\bp}   {\begin{prop}}
\newcommand{\ep}   {\end  {prop}}
\newcommand{\bc}   {\begin{cor}}
\newcommand{\ec}   {\end  {cor}}
\newcommand{\bd}   {\begin{defn}}
\newcommand{\ed}   {\end  {defn}}

\newcommand{\ba}   {\begin{array}}
\newcommand{\ea}   {\end  {array}}
\newcommand{\be}   {\begin{enumerate}}
\newcommand{\ee}   {\end  {enumerate}}
\newcommand{\bi}   {\begin{itemize}}
\newcommand{\ei}   {\end  {itemize}}

\def\eq#1\en{\begin{equation}#1\end{equation}}  
\def\eqsplit#1\ensplit{
	\begin{equation}\begin{split}#1\end{split}\end{equation}
	}
\def\eqalign#1\enalign{
	\begin{align}#1\end{align}
	}
\def\eqmul#1\enmul{
	\begin{multline}#1\end{multline}
	}
\newcommand{\eqarrstar} {\begin{eqnarray*}} 
\newcommand{\enarrstar} {\end{eqnarray*}} 
\newcommand{\eqarray}   {\begin{eqnarray}} 
\newcommand{\enarray}   {\end{eqnarray}}

\newcommand{\lbeq}[1]  {\label{e:#1}}
\newcommand{\refeq}[1] {\eqref{e:#1}}    % AMS-LaTeX trick!
\makeatletter
\newcommand{\labelcounter}[2]{{%
	\stepcounter{#1}%	First, increase the ``countC'' by one.
	\protected@write\@auxout{}%
	{\string\newlabel{#2}{{\csname the#1\endcsname}{\thepage}}}%
	{\ref{#2}}%	Finally, make sure to refer to this label, 
	}}
\makeatother
\newcommand{\sss}   { \scriptscriptstyle } 

\newcommand{\Nbold} {{\mathbb N}}

\newcommand{\Rbold} {{\mathbb R}}

\newcommand{\Zbold} {{\mathbb Z}}

\newcommand{\spose}[1] {{\hbox to 0pt{#1\hss}} }
\newcommand{\ltapprox} {\mathrel{\spose{\lower 3pt\hbox{$\mathchar"218$}}
 \raise 2.0pt\hbox{$\mathchar"13C$}}}
\newcommand{\gtapprox} {\mathrel{\spose{\lower 3pt\hbox{$\mathchar"218$}}
 \raise 2.0pt\hbox{$\mathchar"13E$}}}

\newcommand{\prodtwo}[2]{
	\prod_{ \mbox{ \scriptsize 
		$\begin{array}{c} 
		{#1} \\ 
		{#2}  
		\end{array} $ } 
		} 
	} 
\UseSection  %Necessary to define Numbering scheme for Theorem, etc.
\newcommand{\hlf}{\frac{1}{2}}
\newcommand{\ra}{\rightarrow}

\renewcommand{\to}      {\rightarrow}

\newcommand{\bigcupd} {\stackrel{\cdot}{\bigcup}}
\newcommand{\del}{\partial}

\newcounter{countC}  % Defined the counter ``countC''
\newcounter{countR}  % Defined the counter ``countR''
\newcommand{\R}{\Rbold}

\newcommand{\Z}{\Zbold}
\newcommand{\N}{\Nbold}

\newcommand{\nn}{\nonumber}

\newcommand{\smallsup}[1] {{\scriptscriptstyle{({#1}})}}

\newcommand{\walk}{\vec{x}}
\newcommand{\ewalk}{\vec{\eta}}
\newcommand{\walkvec}[2]{\vec{\eta}^{\smallsup{#1}}_{#2}}
\newcommand{\walkcoor}[2]{\eta^{\smallsup{#1}}_{#2}}

\newcommand{\wt}[1]{\widetilde{#1}}

\newcommand{\mc}[1]{\mathcal{#1}}
\newcommand{\mP}{\mathbb{P}}

\newcommand{\mE}{\mathbb{E}}

\newcommand{\blank}[1]{}
\newcommand{\eps}{\varepsilon}

\newcommand{\Di}{{\cal D}}

\newcommand{\Mi}{{\cal M}}

\renewcommand{\P}{{\mathbb P}}
\newcommand{\E}{{\mathbb E}}

\title  {
        A monotonicity property for random walk\\ in a partially random environment
        }

\author{
Mark Holmes\footnote{Department of Statistics, University of Auckland.  E-mail {\tt
holmes@stat.auckland.ac.nz}} ~and Rongfeng Sun\footnote{Department of Mathematics, National University of Singapore.  E-mail {\tt
matsr@nus.edu.sg}}}

\begin{document}

\maketitle

\begin{abstract}
We prove a law of large numbers for random walks in certain kinds of i.i.d.~random environments in $\Z^d$ that is an extension of a result of Bolthausen, Sznitman and Zeitouni~\cite{BSZ03}.  We use this result, along with the lace expansion for self-interacting random walks, to prove a monotonicity result for the first coordinate of the speed of the random walk under some strong assumptions on the distribution of the environment.


\medskip

\end{abstract}

%\tableofcontents

%\title{A law of large numbers for a class of random walks in random environment}
%\date{}
%\maketitle

\section{Introduction}
\label{sec-intro}
Since the pioneering work of Solomon and others in the mid 1970s to early 1980s, random walk in a random environment (RWRE) has enjoyed a revival in recent times as a number of interesting results have been obtained.  Many of these results relate to laws of large numbers and invariance principles for i.i.d.~random environments that are uniformly elliptic (all nearest-neighbour transition probabilities are bounded away from zero).  While the behavior of one-dimensional RWRE is quite well understood, understanding is much less complete for
RWRE in higher dimensions, and in particular, for non-ballistic RWRE. Several special classes of RWRE are amenable to analysis in general dimensions, such as random walk among random conductances, random walks in balanced random environments or Dirichlet random environments, random walks in random environments which are small perturbations of a deterministic environment, etc (see e.g.~\cite{BarD10, BerD11, ES06, Sab04} and the references therein). Another class is a random walk in a {\em partially random environment}, introduced by Bolthausen, Sznitman and Zeitouni in \cite{BSZ03}. They established laws of large numbers and central limit theorems for a RWRE in dimensions $d=d_0+d_1$, where $d_0\ge 1$ is the number of coordinates in which the environment is random, and where the projection of the walk onto the remaining $d_1\ge 5$ coordinates is a deterministic symmetric random walk.

In this paper we consider monotonicity properties of the speed for random walks in partially random environments (RWpRE) that are similar to those considered in \cite{BSZ03}. Such properties have not been extensively studied in the literature. The few results on the monotonicity of the speed we are aware of include the work of Holmes and Salisbury \cite{HS_RWDRE} where monotonicity of the speed (when it exists) is proved for all environments that take only two possible values via a coupling argument, and Sabot~\cite{Sab04} and Fribergh~\cite{F10}, where asymptotic expansions for the speed are derived for random walks in random environments that are small perturbations of a simple random walk with drift. There has also been recent progress in the study of the monotonicity of the speed as a function of the bias for a biased random walk on supercritical percolation clusters~\cite{FH11} and Galton-Watson trees~\cite{bAFS11}. Our main result is a monotonicity result for the first coordinate of the speed, under some special assumptions on the distribution of the partially random environment.  For example when at each site either the left or right step in the the first coordinate direction is not available, we prove that the first coordinate of the speed is monotone increasing in the probability that the right step is available. Our proof consists of two steps. We first extend a result of \cite{BSZ03} to show that the (non-random) speed exists almost surely for the class of RWpRW under our consideration. We then establish the desired monotonicity by analyzing an expansion formula for the speed derived in~\cite{HH07} using lace expansion techniques, which is valid for all annealed RWRE, but is most useful in the case $d_1\gg d_0$ when one has good control (in terms of finite random walk Green's functions) over the terms in the expansion.

Let $\Mi_1(\Z^d)$ be the space of probability kernels on $\Z^d$, and more generally for $c>0$ let $\Mi_c(\Z^d)$ denote the space of kernels on $\Z^d$ with total mass $c$.  Given a family of probability kernels $\omega:=(\omega_{x,m}(\cdot))_{x\in\Z^d, m\in\N} \in\Mi_1(\Z^d)^{\Z^d\times\N}$ which we call a {\em cookie environment}, the law of a random walk $(X_n)_{n\geq 0}$ in the cookie environment $\omega$ starting at $X_0=x$, denoted by $P_{x,\omega}$, is defined as follows.  Under $P_{x,\omega}$, the walk evolves conditional on its history via the transition probabilities
$P_{x,\omega}(X_n=X_{n-1}+u |(X_i)_{0\leq i\leq n-1}) = \omega_{X_{n-1}, \ell_{n-1}(X_{n-1})}(u)$,
for all $n \in\N$, where $\ell_n(y)= \sum_{k=0}^n 1_{\{X_k=y\}}$ is the number of visits to $y$ up to time $n$. In words, upon the $m$-th visit to $x$, the walk
sees the environment $\omega_{x,m}$ and makes a jump accordingly. We will consider the case when $\omega$ is random, and the cookie environment at different points in space, $(\omega_{x,\cdot}(\cdot))_{x\in\Z^d}$, are i.i.d.\ with a common law $\mu\in \Mi_1(\Mi_1(\Z^d)^{\N})$.  The measure $P_{x,\omega}$ is called the {\em quenched} law.  When we average the quenched law of  $(X_n)_{n\geq 0}$ with respect to the cookie environment $\omega$, we obtain the so-called {\em annealed}
(or more accurately the {\em averaged}) law
$$
P_x := \P\times P_{x,\omega},
$$
where $\P:=\mu^{\otimes \Z^d}$ denotes the law of $\omega$ in the product space $\Omega:=(\Mi_1(\Z^d))^{\Z^d\times\N}$.

The random walk model described above is sometimes called a multi-excited random walk (in a random cookie environment).  When we restrict the cookie environment to environments that are constant in $m$, i.e.~$\omega_{x,m}(\cdot)\equiv \omega_{x}(\cdot)$ for every $x$ and $m$, we obtain the more often studied RWRE model.  Among other assumptions required for our two main results, the result on monotonicity of the speed assumes that $\omega_{x,m}(\cdot)\equiv \omega_{x}(\cdot)$, while the law of large numbers result does not.

\subsection{The law of large numbers}
\label{sec:LLN}
Recall that the RWRE $X$ is said to satisfy a (strong) law of large numbers (LLN), if there exists a constant $v\in \R^d$, such that $P_o$ a.s.\ $\lim_{n\to\infty} \frac{X_n}{n}=v$. For ballistic RWRE, the technique of regeneration times has been useful in proving the LLN, see e.g.~\cite{SZ99, Sz04, Zeit04}.
For non-ballistic RWRE, it is not known in general if there is a deterministic limit, however Sznitman and Zerner~\cite{SZ99} and Zerner~\cite{Zern02} showed that $X_n/n$ converges $P_o$-almost surely to a random variable taking at most two possible values. In dimension two, the LLN has been established by Zerner and Merkl~\cite{MZ01}.  In both cases the environment is assumed to be i.i.d.~and elliptic (with $\omega_{x,m}(\cdot)\equiv\omega_{x}(\cdot)$), but it was shown recently by Holmes and Salisbury in~\cite[Sec.~3]{HS_RWDRE} that the assumption of ellipticity can be dropped. In higher dimensions, the LLN has only been established for various special classes which do not require ballisticity. One such special class is a random walk in a partially random environment studied by Bolthausen, Sznitman, and Zeitouni in \cite{BSZ03}. They assume that $d=d_0+d_1$ with $d_0\geq 1$ and $d_1\geq 5$, and there exists $\alpha \in (0,1)$ and $q(\cdot) \in \Mi_1(\Z^{d_1})$, such that for $\mu$ a.e.\ $\omega_x(\cdot)\in\Mi_1(\Z^d)$:
\begin{itemize}
\item[\rm (a)] $\omega_x(\cdot)$ is supported on the canonical unit vectors, i.e., $\sum_{e\in\Z^d, |e|=1} \omega_x(e)=1$,

\item[\rm (b)] if $\Pi_{d_1}: \Z^{d}\to \Z^{d_1}$ denotes the projection of $v\in\Z^d$ to its last $d_1$ coordinates, then
\eq\lbeq{alpha}
(\omega_x\circ \Pi_{d_1}^{-1})(\cdot) = \alpha q(\cdot) + (1-\alpha)\delta_0(\cdot),
\en
\item[\rm (c)] $\sum_{e\in \Z^{d_1}, |e|=1} q(e)=1$, and $q(e)=q(-e)>0$ for all $e\in\Z^{d_1}$ with $|e|=1$.
\end{itemize}
When the RWRE $(X_n)_{n\geq 0}$ on $\Z^d$ is projected to its last $d_1$ coordinates, one obtains a random walk
$(Y_n)_{n\in\N}:=(\Pi_{d_1} X_n)_{n\in\N}$
on $\Z^{d_1}$ with transition kernel $\alpha q(\cdot) + (1-\alpha)\delta_0(\cdot)$. In dimensions $d_1\geq 5$, such a random walk admits so-called {\em cut times},
i.e.,
\eq
\lbeq{cuttimes}
\Di := \{ n\in \Z : Y_{(-\infty, n-1]} \cap Y_{[n,\infty)} = \emptyset\} \neq \emptyset \qquad a.s.
\en
Using the fact that conditional on the projected random walk $(\Pi_{d_1} X_n)_{n\geq 0}$, the environments the walk $X$ sees between successive cut times are independent,
a law of large numbers was proved in \cite{BSZ03}. The key assumption is thus the existence of cut times, which play the role of regeneration times in this context.

We now extend the aforementioned LLN from~\cite{BSZ03} to cookie environments $\omega:=(\omega_{x,m}(\cdot))_{x\in\Z^d, m\in\N}$ where $(\omega_{x,\cdot}(\cdot))_{x\in\Z^d}$ are i.i.d.\ with common law $\mu$. Furthermore, we will allow $\alpha$ in \refeq{alpha} to be random. More precisely,
conditions (a)--(c) will be relaxed to the following:
\begin{itemize}
\item[\rm (a')] There exists $K>0$ such that $\omega_{x,m}(\cdot)$ is a.s.\ supported on $[-K,K]^d$ for all $m\in\N$.

\item[\rm (b')] There exists some $\delta>0$ and $q\in \Mi_1(\Z^{d_1})$ with $q(0)=0$ such that, for all $n\in\N$, a.s.\
\eq\lbeq{muproject}
(\omega_{x,m}\circ \Pi_{d_1}^{-1})(\cdot) = \alpha_{\omega_{x,m}} q(\cdot) + (1-\alpha_{\omega_{x,m}})\delta_0(\cdot) \qquad \mbox{ for some } \alpha_{\omega_{x,m}}\geq \delta.
\en

\item[\rm (c')] For a random walk $(Y_n)_{n\in\N}$ on $\Z^{d_1}$ with jump probability kernel $q(\cdot)$, $\Di\neq\emptyset$ a.s.
\end{itemize}
Note that condition (c') holds for any $Y$ with a non-zero drift, and it also holds for any $Y$ that is at least 5-dimensional in the sense that the space spanned by vectors in the support of $Y$'s increment distribution is at least 5-dimensional (see e.g.\ \cite[Lemma 1.1]{BSZ03}). By ergodicity, (c') in fact implies that $\cal D$
is an infinite set almost surely.

\begin{THM}
\label{T:LLN}{\bf (Law of large numbers)} Let $(X_n)_{n\geq 0}$ be a random walk in a cookie environment $\omega$, where $(\omega_{x,\cdot}(\cdot))_{x\in\Z^d}$ are i.i.d.\ with common law $\mu$ and satisfy conditions {\rm (a')}--{\rm (c')} above. Then there exists some $v\in \R^d$ such that
\eq\lbeq{LLN}
\lim_{n\to\infty} \frac{X_n}{n} = v \qquad P_o \mbox{ almost surely}.
\en
\end{THM}
%\begin{REM}
%It will be clear from the proof that Theorem~\ref{T:LLN} also applies to random environments where given an environment $\omega$ as in %Theorem~\ref{T:LLN},
%the transition probabilities of the walk are given by $\P(X_n=X_{n-1}+v |(X_i)_{0\leq i\leq n-1}) = \omega_{X_{n-1}, l_{n-1}(X_{n-1})}(v)$, where %$l_k(x)=\sum_{i=0}^k 1_{\{X_i=x\}}$ is the local time of $X$ at $x$ at time $k$. For such a walk, the transition kernel $\omega_{x,k}$ is used upon %the $k$-th visit to $x$. Multi-excited
%random walk is one such example.  This version of the extension is referred to and applied in \cite{H08compete}.
%\end{REM}
Although we will only need time-independent random environments when we later study the monotonicity of the speed $v$, we have formulated Theorem~\ref{T:LLN}
for cookie environments because there are natural interesting examples involving cookie environments, such as the (multi)-excited random walk (see e.g.~\cite{H08compete}). Apart
from extending \cite[Theorem 1.4]{BSZ03}, Theorem \ref{T:LLN} also extends the weak law of large numbers \cite[Theorem 2.5]{HH07}, which incorporated the parameter $\delta$ but is only valid in the perturbative regime where there is some $\epsilon>0$ sufficiently small (depending on $\delta$) such that $\mP(|\omega_o-\mE[\omega_o]|_1>\epsilon)=0$, where $|\cdot|_1$ denotes here the total variation norm on $\Mi_1(\Z^d)$.

The proof of Theorem~\ref{T:LLN} is based on adaptations of arguments in \cite{BSZ03}, which we outline in Section \ref{sec:LLNproof}. Relaxing conditions (a) and (c) to
(a') and (c') does not induce any change in the proof. However relaxing condition (b) to (b') requires a simple but not entirely trivial observation, and the proof needs to be modified accordingly. Indeed we were encouraged to consider this extension by one of the authors of \cite{BSZ03}.

The lace expansion for self-interacting random walks of van der Hofstad and Holmes \cite{HH07} gives the following series representation for the expected increment of the RWRE under $P_o$,
\eqalign
E_o[X_n-X_{n-1}]=E_o[X_1] +\sum_{m=2}^n \sum_x x\pi_m(x)\lbeq{incformula},
    \enalign
 where $\pi_m(x)$, for $m\ge 2, x\in \Z^d$ are somewhat complicated quantities known as lace-expansion coefficients.   If \refeq{incformula} converges then $E_o[X_n/n]$ converges to the same limit (convergence of $E_o[X_n/n]$ also follows from dominated convergence and the fact that $X_n/n$ converges almost surely).
Theorem \ref{T:LLN} allows us to go one step further and say that if \refeq{incformula} converges then it converges to $v$, i.e.
 \eqalign
    v& = E_o[X_1] +\lim_{n\ra \infty}\sum_{m=2}^n \sum_x x\pi_m(x)\lbeq{vformula}.
    \enalign
Usually in analysing this formula we first require enough control on the coefficients $\pi_m(x)$ to ensure that this series converges.
%Once Theorem \ref{T:LLN} has been proved, the lace expansion for self-interacting random walks of van der %Hofstad and Holmes \cite{HH07} gives an infinite series representation for the velocity $v$, provided that the %series converges.  By Proposition 3.1 of \cite{HH07}, the velocity is given by
%    \eqalign
%    v& = \lim_{n\ra \infty}E_o[X_n-X_{n-1}]=E_o[X_1] +\lim_{n\ra \infty}\sum_{m=2}^n \sum_x %x\pi_m(x)\lbeq{vformula},
%    \enalign
% where $\pi_m(x)$, for $m\ge 2, x\in \Z^d$ are somewhat complicated quantities known as lace-expansion %coefficients.
To prove other properties of the velocity such as the sign, continuity, differentiability, and monotonicity, we typically require even better control of the coefficients $\pi_m(x)$.  The expansion is not perturbative in the usual sense.  For any translation invariant self-interacting random walk (see \cite{HH07} for precise details) for which $E_o[X_n-X_{n-1}]$ converges, the formula \refeq{incformula} is valid.  Currently in order to extract useful information from the formula we require the walk to be sufficiently and quantifiably transient, independent of the history of the walk.  Analysis of the formula without such an assumption would require a major advance in our understanding of the expansion methodology itself.  In this paper the walks that we consider have a high-dimensional simple random walk component (see e.g.~assumption {\bf (A3)} in Section \ref{sec-mono}) which has the necessary transience property.

We will study \refeq{vformula} as a function of a particular parameter of interest, $\beta$, under some strong assumptions on the distribution of the environment, and show that the velocity $v^{[1]}$ (the first coordinate of the velocity) is increasing in $\beta$.  A precise formulation of this result is given in Theorem~\ref{thm:monotone} at the end of Section \ref{sec-mono}.

\subsection{Monotonicity}
\label{sec-mono}
Our discussion and results concerning monotonicity are restricted to time-independent random environments, i.e., $\omega_{x,m}(\cdot)\equiv \omega_x(\cdot)$
for all $m\in\N$ and $x\in\Z^d$.  For a discussion of monotonicity in cookie environments with $d_0=1$ and condition (b) instead of (b'), see for example \cite{H08compete}.  Rather a lot is known when $d=1$, see Holmes and Salisbury \cite{HS_comb}.

It is well known that the velocity of a RWRE is not monotone increasing in the expected local drift at the origin.  In fact it is possible in one dimension that the expected local drift is negative, yet the walk is transient to the right.  In higher dimensions (see e.g.~\cite{BSZ03},\cite{HS_RWDRE}) the speed and the expected local drift can even carry opposite signs. For example, consider a nearest-neighbour RWRE on $\Z^2$ with i.i.d.\ environment $\omega:=(\omega_x)_{x\in\Z}$, where $\P(\omega_o(e_1)=\omega_o(e_2)=\hlf)=p$ and $\P(\omega_o(-e_1)=1)=1-p$.  Due to an elementary renewal structure, the velocity of this random walk can be calculated explicitly \cite{HS_RWDRE} as
\[v=\left(\frac{p(2-p)}{2+3p-2p^2-p^3}\right)\cdot (3,1)-(1,0).\]
It is easy to see that the first coordinate $v^{[1]}$ need not carry the same sign as the expected local drift $\frac{3p}{2}-1$.  However in this example the velocity is indeed monotone in $p$.  Holmes and Salisbury \cite{HS_RWDRE} prove that this is the case for any 2-valued environment.  To be precise, if $\mP(\omega_o=A_1)=p=1-\mP(\omega_o=A_2)$ then the velocity $v(p)$ if it exists is monotone in $p$.  This fails in general \cite{HS_RWDRE} for 3-valued environments with respective probabilities $p(1-q)$, $(1-p)(1-q)$ and $q$ for fixed $q$.

%in one dimension, consider a nearest-neighbour RWRE on $\Z$ with i.i.d.\ environment $\omega:=(\omega_x)_{x\in\Z}$, where $\P(\omega_o(e_1)=1)=.25=1-\P(\omega_o(e_1)=.25)$.  Then the (annealed) expected local drift is $E_o[X_1]=.25+.25\times .75-.75^2<0$, but the positive density of sites at which it is impossible to go left ensures that the walk is transient to the right. Because $\E[\omega_o(-e_1)/\omega_o(e_1)]^{-1}<1<\E[\omega_o(e_1)/\omega_o(-e_1)]$, the speed of this walk is almost surely zero  by the classic result of Solomon (see e.g.~\cite{Sz04}).

Now consider two i.i.d.\ environments $\omega$ and $\wt{\omega}$ satisfying conditions (a)--(c) in Sec.~\ref{sec:LLN} with the same $d_0, d_1, \alpha$ and $q(\cdot)$.
Assume furthermore that for some $\kappa\in (0,1)$,
\begin{equation}\label{kap}
\omega_o(e_1)+\omega_o(-e_1)=\wt{\omega}_o(e_1)+\wt{\omega}_o(-e_1)=\kappa \qquad \P \ a.s.,
\end{equation}
and $\wt{\omega}_o(e_1)$ stochastically dominates $\omega_o(e_1)$ in the sense that
\begin{equation}\label{dom}
\P(\omega_o(e_1)\le s)\geq \P(\wt{\omega}_o(e_1)\le s)  \qquad \mbox{for all } s\in [0,1].
\end{equation}
Suppose that a random walk $X$ (resp.\ $\wt{X}$) with $X_0=\wt{X}_0=0$ in the random environment $\omega$ (resp.\ $\wt{\omega}$) satisfies a LLN with (deterministic) speed $v$ (resp.\ $\wt{v}$), is it true that $v^{[1]} \leq \wt{v}^{[1]}$, where $w^{[1]}$ denotes the first coordinate of $w\in \R^d$?

When $d_0=1$ and $d_1\ge 0$, the answer to the above question is affirmative, since we can easily couple $(X, \omega)$ with $(\wt{X},\wt{\omega})$ such that at each time $n$, $\wt{X}_n^{[i]}=X_n^{[i]}$ for $i>1$, and $\wt{X}_n^{[1]}\ge X_n^{[1]}$.  That is, the position of the walks $\{X_n\}_{n\ge 0}$ and $\{\wt{X}_n\}_{n\ge 0}$ differ only in the first coordinate and $\wt{X}$ is never to the left of $X$.
%The answer is also affirmative when the environment takes only two possible values \cite{HS_RWDRE}.

In general, however, we expect the answer to the above question to be negative, since the limiting velocity $v$ depends on the joint distribution
of $(\omega_o(\pm e_i))_{1\leq i\leq d_0}$. Consider for example the case $d_0=d=2$ with $\kappa=\hlf$, such that
\begin{equation}\label{ex1}
\P\big((\omega_o(e_i),\omega_o(-e_i))=(1/2,0)\big)=\beta=1-\P\big((\omega_o(e_i),\omega_o(-e_i))=(0,1/2)\big), \qquad i=1,2,
\end{equation}
for some $\beta\in (0,1)$, and $(\omega_o(e_2),\omega_o(-e_2))$ is independent of $(\omega_o(e_1),\omega_o(-e_1))$.  The corresponding random walk has limiting velocity
$v=0$ for all $\beta\in (0,1)$, since it eventually gets stuck on a finite number of sites (e.g.~see \cite{HS_RWDRE}).  On the other hand one expects that the random walk in the random environment $\wt{\omega}$ (again with $d_0=d=2$, $\kappa=\hlf$) such that
\begin{equation}\label{ex2}
\P\big(\wt{\omega}_o(e_2)=\wt{\omega}_o(e_1)=1/2\big)=\beta=1-\P\big(\wt{\omega}_o(-e_2)=\wt{\omega}_o(-e_1)=1/2\big)
\end{equation}
has a non-trivial deterministic velocity
%in each coordinate direction
whenever $\beta\ne 1/2$ (see e.g.\ \cite{HS_RWDRE}).

%It is also worth noting that there may be parametric families of environments where monotonicity holds, even though there is no stochastic domination.  For example, let $d_0=1$, $c_1=d^{-1}$, $\beta \in [0,1]$ and $\mu(X=\frac{1+\beta}{2})=\frac{2+\beta}{3}=1-\mu(X=\frac{1-\beta}{2})$.  We expect that the speed of a random walk in this environment is monotone increasing in $\beta$ when the dimension $d$ is sufficiently large.

%In the presence of stochastic domination in the first coordinate, does this monotonicity property hold when the probability of stepping in the first coordinate is not constant?  We expect not.

As a special case of our main result (Theorem~\ref{thm:monotone}) below, we will show that for each of the two examples in (\ref{ex1}) and (\ref{ex2}), if $ d_1\gg 1$ and
the $1/2$ in (\ref{ex1}) and (\ref{ex2}) is replaced by a sufficiently small constant, then $v^{[1]}$ is continuous and strictly increasing in $\beta$.  We believe that one can tune parameters in these two examples such that $\wt{\omega}_o(e_1)$ stochastically dominates $\omega_o(e_1)$ as in (\ref{dom}), and yet $v^{[1]}>\wt{v}^{[1]}$.
%RWRE where $\omega_o$ is $\P$-a.s.\ the uniform distribution on a random subset of $\{\pm e_i\}_{1\leq i\leq d}$ is studied more extensively in \cite{HS_RWDRE}.

%If one then assumes that for fixed $\beta$, the asymptotic velocities are different (say $v_{X,\wt{X}}>v_{X,X}$), then from the continuity of the velocities as functions of $\beta$ we see that for some $\beta'<\beta$, $v_{X',\wt{X}}>v_{X,X}$, where $X'$ is independent of $\wt{X}$ and has Bernoulli distribution with parameter $\beta'$.

For $\omega$ and $\wt{\omega}$ formulated as in (\ref{kap}) and (\ref{dom}), if we further assume that $\{\omega_o(\pm e_1)\}$ is independent of $\{\omega_o(\pm e_i) :2\leq i\leq d\}$, the same holds for $\wt{\omega}$, and $\{\omega_o(\pm e_i) :2\leq i\leq d\}$ is equally distributed with $\{\wt{\omega}_o(\pm e_i) :2\leq i\leq d\}$, then it is natural to expect the monotonicity property $v^{[1]}\leq \wt{v}^{[1]}$ to hold. We will prove this in high dimensions in some special cases.
\bigskip

We now formulate precisely the class of RWRE for which we can prove monotonicity of $v^{[1]}$. Let $d=d_0+d_1$, $1\leq d_*\leq d_0$, and let $\gamma,\kappa,\delta\in (0,1]$ be such that $\gamma +\delta\le 1$. Our assumptions on the random environment $\omega$ consist of the following.
\medskip

\noindent{\bf (A0):} $\omega:=(\omega_x)_{x\in\Z^d}$ are i.i.d.\ with common law $\mu\in \Mi_1(\Mi_1(\Z^d))$, and $\mu$-a.s., $\omega_o$ is supported
on $V_d:=\{\pm e_i\}_{1\leq i\leq d}$ and satisfies conditions (b')--(c') in Sec.~\ref{sec:LLN} for some $q(\cdot)\in\Mi_1(\Z^{d_1})$.
\medskip

\noindent{\bf (A1):} $\omega_o$ as an $\Mi_1(\Z^d)$-valued random variable furthermore admits the representation
\begin{equation}\label{omegarep}
\omega_o= \xi_{d_*}\times \delta_{d-d_*}(o) + \delta_{d_*}(o)\times \tilde \xi_{d-d_*},
\end{equation}
where $\times$ denotes product measure, $\delta_i(o)\in \Mi_1(\Z^i)$ denotes the delta measure at the origin $o$, and $\xi_{d_*}$ and $\tilde \xi_{d-d_*}$ are {\em independent} random kernels supported on $V_{d_*}$ (resp.\ $V_{d-d_*}$) with total mass $\gamma$ (resp.\ $1-\gamma$) and laws $\mu_{d_*}\in \Mi_1(\Mi_\gamma(\Z^{d_*}))$ (resp.\ $\tilde \mu_{d-d_*}\in \Mi_1(\Mi_{1-\gamma}(\Z^{d-d_*}))$).
\medskip

\noindent{\bf (A2):} There exist $\nu_1, \nu_2 \in \Mi_\gamma(\Z^{d_*})$ with {\em disjoint supports} $\mc{S}_1, \mc{S}_2\subset V_{d_*}$ such that
\begin{equation}\label{mud*}
\begin{aligned}
&\mu_{d_*}(\xi_{d_*}=\nu_1) =  \kappa(1-\beta), \qquad \mu_{d_*}(\xi_{d_*}=\nu_2) = \kappa \beta, \\
&\mu_{d_*}\big({\rm supp}(\xi_{d_*})\cap \big(\mc{S}_1\cup \mc{S}_2\cup \{\pm e_1\}\big)=\emptyset\big) =  1-\kappa,
\end{aligned}
\end{equation}
and
\eq
\lbeq{rhodef}
\rho:=(\nu_2(e_1)-\nu_2(-e_1))-(\nu_1(e_1)-\nu_1(-e_1))>0.
\en
For simplicity, we will also use $\mc{S}_i$ to denote $\mc{S}_i\times 0_{d-d_*}\subset\Z^d$, where $0_k$ is the zero vector in $\Z^k$.
\medskip

\noindent{\bf (A3):} Let $q(\cdot)\in\Mi_1(\Z^{d_1})$ be as in (b'), and let $G_q(x):=\sum_{k=0}^\infty q^{*k}(x)$, where $q^{*k}(\cdot)$ denotes the $k$-fold
convolution of $q(\cdot)$ with itself. Then
\begin{equation}\label{Gcond}
G_q(o)<2 \qquad \mbox{ and } \qquad G^{*i}_q:=\sup_{x\in\Z^{d_1}} G_q^{*i}(x)<\infty \mbox{  for  } i=1,2,3,4.
\end{equation}

Let us elaborate more on the assumptions. Assumption {\bf (A1)} requires that $\mu$-a.s., a walk with transition kernel $\omega_o$ will make a nearest-neighbor jump
in the first $d_*$ coordinate directions with probability $\gamma$, and make a nearest-neighbor jump in the last $d-d_*$ coordinate directions with probability $1-\gamma$.
Furthermore, the restriction of $\omega_o$ to unit vectors in the first $d_*$ coordinate directions is independent of its restriction to unit vectors in the last
$d-d_*$ coordinate directions.

The parameter $\beta$ in {\bf (A2)} allows us to tune the relative weight of the kernels $\nu_1$ and $\nu_2$. Since $\rho>0$, $\nu_1$ and $\nu_2$ have different drifts in the first coordinate direction. We can therefore expect $v^{[1]}$ to increase monotonically as we increase the weight of $\nu_2$ at the expense of $\nu_1$.
However, our proof requires $\nu_1$ and $\nu_2$ to have disjoint supports. Furthermore, when $\xi_{d_*}$, the restriction of $\omega_x$ to $V_{d_*}$, is neither $\nu_1$ nor $\nu_2$, then its support must be disjoint from the support of $\nu_1$ and $\nu_2$, and it cannot allow jumps that change the first coordinate
of the random walk (which is trivially satisfied if $\kappa=1$). Under this assumption, the history of the walk up to time $n$  either provides no information about $\xi_{d_*}$ at $x\in\Z^d$ because no jumps affecting the first $d_*$ coordinates has been taken from $x$, or we can determine whether $\xi_{d_*}=\nu_1$, or $\nu_2$, or neither, based on past jumps from $x$ affecting the first $d_*$ coordinates. This fact will be crucial for our proofs.

Assumption {\bf (A0)}, and in particular (b'), guarantees that
for $\mu$ almost every realization of $\omega_o$, with probability at least $\delta$, a jump following $\omega_o$ will induce a change in (and only in) the last $d_1$
coordinates, and conditional on this event, the jump follows the kernel $q$. This allows us to extract
a deterministic random walk and apply lace expansion techniques.

Assumption {\bf (A3)} is needed to control the lace expansion coefficients. It is true for example when $d_1$ is sufficiently large.  As an alternative to the assumption $G^{*4}<\infty$ in {\bf (A3)}, we can instead assume that an appropriate local central limit theorem type bound, \refeq{localCLT}, holds for the random walk with kernel $q$.  See Section \ref{sec:G4} for further details. We note that when $q(\cdot)$ has zero mean, $G^{*4}<\infty$ if and only if $d_1>8$, while the local limit theorem bound \refeq{localCLT} holds as long as $d_1>6$. We expect that the methods of this paper could be adapted to handle cases where $d_1$ is small for {\em asymmetric} $q$, provided that the bias of $q(\cdot)$ is sufficiently strong.  This analysis would require different estimates, similar to those used in the analysis of once-reinforced random walk with drift in \cite{HH07}.
%(when e.g.~$d_1=1$, {\bf (A3)} is not satisfied even if $q$ is asymmetric)
%If we replace condition (b') with (b) then we do not require the assumption $G^{*4}<\infty$ in {\bf (A3)}.
\bigskip

We are now ready to state the second main result of this paper, which holds under the further assumption that $\delta$ in \refeq{muproject} is sufficiently
close to 1.

\begin{THM}{\bf (Monotonicity of speed)}
\label{thm:monotone}
Let $X$ be a random walk in an i.i.d.\ random environment $\omega$ which satisfies {\bf (A0)--(A4)} with kernel $q(\cdot)$ and constants $\gamma, \kappa, \delta \in (0,1]$, $\gamma+\delta\leq 1$. There exists $\delta_q\in (0,1)$ depending only on $q(\cdot)$, such that if $\delta\in (\delta_q,1)$, then $v^{[1]}$ is continuous and strictly increasing in $\beta$.
\end{THM}
Note that under assumptions {\bf (A0)--(A4)}, the existence of a deterministic velocity $v$ is guaranteed by Theorem~\ref{T:LLN}.

The simplest random environment for which Theorem~\ref{thm:monotone} applies is when: (1) $d_*=d_0$, which implies that in (\ref{omegarep}), $\tilde \xi_{d-d_*}=\tilde \xi_{d_1}=q(\cdot)$ a.s., and hence that the probability of taking a $q$-step (a step with non-zero $\Pi_{d_1}$ projection) is the constant $\delta=1-\gamma$; (2) $\kappa=1$ so that $\xi_{d_*}$ in (\ref{omegarep}) equals either $\nu_1$ or $\nu_2$. Such a random environment allows only two possible realizations for each $\omega_x$, and $\beta$ determines the probabilities of their occurrence, so we know from \cite{HS_RWDRE} that the velocity is monotone in $\beta$.  The more general random environments formulated in Theorem~\ref{thm:monotone} can be regarded as perturbations of this simple case by allowing more randomness: $\kappa<1$ allows $\xi_{d_*}$ to take on realizations other than $\nu_1, \nu_2$, provided assumption {\bf (A2)} holds; $d_*<d_0$ allows $\tilde \xi_{d-d_*}$ to be random, provided condition (b') holds.
We note that suitable assumptions on the additional randomness is necessary for the monotonicity result to hold, as it was shown in~\cite{HS_RWDRE} that monotonicity does not hold in general for random environments which almost surely takes on one of three possible realizations. 

Lastly we remark that in related works of Sabot~\cite{Sab04} and Fribergh~\cite{F10}, the authors study the speed of a RWRE which is a perturbation of a homogeneous simple random walk with drift. This is similar in spirit to our model since both are perturbations of a simple random walk. However the exact nature of perturbations, the results, and the techniques are all quite different. In~\cite{Sab04} and~\cite{F10}, the authors rely on the representation of the speed of a ballistic random walk in terms of Kalikow's auxiliary random walk and perform perturbation expansion for its Green function. This results in an expansion of the speed as a function of the perturbation parameter, which implies monotonicity of the speed when the perturbation parameter is zero. Our results employ lace expansion techniques and require stronger assumptions, but are valid as long as the deterministic component of the random walk is transient enough (not necessarily ballistic), as characterized by the condition on its Green function in {\bf (A3)}. The monotonicity we obtain is valid on the whole range of admissible parameters $\beta$.

\subsection{Organisation}
\label{sec:org}
The remainder of the paper is organised as follows.  In Section \ref{sec:LLNproof} we prove Theorem \ref{T:LLN}.  In Section \ref{sec:lace-exp} we review the relevant notation and results on the lace expansion for self-interacting random walks from \cite{HH07}, including the formula for the speed.  The basic ingredients of the formula are annealed transition probabilities, and these are examined in Section \ref{sec-transprob}.  Sections \ref{sec:formula_converges} and \ref{sec-derivpi} are devoted to proving bounds on some of the quantities appearing in the speed formula and its derivative.  Finally in Section \ref{sec:mono_proof} we prove Theorem~\ref{thm:monotone}
as a consequence of the given formula for the speed and its derivative.

\section{Proof of Theorem \ref{T:LLN}}
\label{sec:LLNproof}

The proof is based on adaptations of the arguments in the proof of \cite[Theorem~1.4]{BSZ03}. The first step is to give a suitable representation of the random cookie environment $\omega$ in terms of independent environments using cut times. Then the standard LLN for i.i.d.\ random variables can be applied.

Let $(Z_i)_{i\in\Z}$ be i.i.d.\ $\Z^{d_1}$-valued random variables with common distribution $q(\cdot)$. Then we can construct a doubly infinite random walk path $(Y_n)_{n\in\Z}$ in $\Z^{d_1}$
with $Y_0=0$, $Y_n = \sum_{i=1}^n Z_i$ for $n\geq 1$, and $Y_n= -\sum_{i=n+1}^0 Z_i$ for $n\leq -1$. The set of {\em cut times} of $(Y_n)_{n\in\Z}$, or
$(Z_n)_{n\in\Z}$, defined in \refeq{cuttimes}, is almost surely non-empty by assumption (c'). Note that $(Z_i, 1_{\{i\in\Di\}})_{i\in\Z}$ is an ergodic sequence with respect to the time shifts
\eq\lbeq{thetak}
\theta_k(Z_i, 1_{\{i\in\Di\}}) = (Z_{i+k}, 1_{\{i+k\in\Di\}}), \qquad \qquad k \in \Z.
\en
Therefore almost surely, $\sup\Di=\infty$ and $\inf\Di=-\infty$. We will denote $\Di\cap [1,\infty) :=\{T_1 < T_2<\cdots\}$ and $\Di\cap (-\infty, 0]:=\{\cdots < T_{-1}<T_0\}$.

We can couple the random walk $X$ in the random cookie environment with $Y$ as follows. Given $(X_n)_{n\geq 0}$, let $\tau_0=0$, and inductively, define
\eq\lbeq{tau}
\tau_j := \inf\{n > \tau_{j-1} : \Pi_{d_1}(X_n-X_{n-1}) \neq 0\}, \qquad j\in\N.
\en
By condition (b') on the random cookie environment $\omega$, $(X_0, X_{\tau_1}, X_{\tau_2},\cdots)$ is distributed exactly as $(Y_n)_{n\geq 0}$. So without loss of generality,
we will assume that $\Pi_{d_1}(X_{\tau_i}-X_{\tau_{i-1}})=Z_i$ for $i\in\N$. Note that when condition (b) holds, $(\Pi_{d_1}X_n)_{n\geq 0}$ is distributed as a random walk $\tilde Y$ on $\Z^{d_1}$ with increment distribution $\alpha q(\cdot) + (1-\alpha)\delta_0(\cdot)$. We can then just couple $X$ with $\tilde Y$ and there will be no
need to introduce the stopping times $\tau_i$. This was the approach taken in \cite{BSZ03}.

By the definition of cut times of $(Z_n)_{n\in\Z}$ and the assumption that $(\omega_{x,m})_{m\in\N}$ is i.i.d.\ for $x\in\Z^d$, we observe that conditional on $(\Pi_{d_1}(X_{\tau_i}-X_{\tau_{i-1}}))_{i\geq 0} = (Z_i)_{i\geq 0}$, the random walk sees spatially distinct, and hence independent, cookie environments on
the time intervals $[0, \tau_{T_1}-1]$, and $[\tau_{T_i}, \tau_{T_{i+1}}-1]$ for $i\in\N$. We thus have the following construction of the annealed law of $X$.

Let $(Z_i)_{i\in\Z}$, $(Y_i)_{i\in\Z}$, and $\Di=(T_i)_{i\in\Z}$ be as above. Let $\omega^{\langle i\rangle}$, $i\in\Z$, be i.i.d.\ $\Omega$-valued random variables
equally distributed with $\omega$, which will be the cookie environment the walk $X$ sees on the time interval $[\tau_{T_{i-1}}, \tau_{T_i}-1]$. For each $i\in\Z$, we
construct a random walk $(X^{\langle i\rangle}_n)_{n\geq 0}$ in the cookie environment $\omega^{\langle i\rangle}$ inductively as follows. Set
$X^{\langle i\rangle}_0=0$. Let
\eq\lbeq{stepnumber}
N(X^{\langle i\rangle}_{[0,n]}) := |\{1\leq j\leq n: \Pi_{d_1}(X^{\langle i\rangle}_j-X^{\langle i\rangle}_{j-1})\neq 0\}|
\en
be the number of steps $X^{\langle i\rangle}$ has taken with non-zero increments in the last $d_1$ coordinates up to time $n$. For each $v\in\Z^d$, if $\Pi_{d_1}v = 0$, then we set $X^{\langle i\rangle}_{n+1}=X^{\langle i\rangle}_n + v$ with probability $\omega^{\langle i\rangle}_{X^{\langle i\rangle}_n,\ell_n(X_n^{\langle i\rangle})}(v)$ so that the walk's jump is unperturbed if it does not affect the last $d_1$ coordinates.
On the event that the walk's jump does affect the last
$d_1$ coordinates, we will change the law of the jump so that its $\Pi_{d_1}$ projection agrees with the next increment of $Y$. More precisely,
for $v$ with
$\Pi_{d_1}v = Z_{T_{i-1}+N(X^{\langle i\rangle}_{[0,n]})+1}$, we assign $v$ as the next jump with probability $\omega^{\langle i\rangle}_{X^{\langle i\rangle}_n,\ell_n(X_n^{\langle i\rangle})}(v)/q(\Pi_{d_1}v)$, which by (b') is equal to
\[\sum_{u:\Pi_{d_1}(u)\ne 0}\omega^{\langle i\rangle}_{X^{\langle i\rangle}_n,\ell_n(X_n^{\langle i\rangle})}(u) \times \frac{\omega^{\langle i\rangle}_{X^{\langle i\rangle}_n,\ell_n(X_n^{\langle i\rangle})}(v)}{\underset{\scriptscriptstyle u:\Pi_{d_1}(u)=\Pi_{d_1}(v)}{\sum}\omega^{\langle i\rangle}_{X^{\langle i\rangle}_n,\ell_n(X_n^{\langle i\rangle})}(u)},\]
i.e.~the product of the probability that there is a jump affecting the last $d_1$-coordinates
and the probability that the jump equals $v$ conditional on $\Pi_{d_1}v = Z_{T_{i-1}+N(X^{\langle i\rangle}_{[0,n]})+1}$. All other jumps are assigned probability $0$. This then couples $X^{\langle i\rangle}$ and $Y$ so that the increments of $X^{\langle i\rangle}$ in the last $d_1$ coordinates match exactly with $(Z_{T_{i-1}+n})_{n\geq 1}$.
% Using an expression of the form
% $$\sum_{v:\Pi_{d_1}(v)\ne 0}\omega(v)=\sum_{z\ne 0}q(z)\sum_{v:\Pi_{d_1}(v)=z}\frac{\omega(v)}{q(\Pi_{d_1}(v))},$$
% we see that the sum of the above probabilities over $v$ add up to $1$.

Since 0 need not be a cut-time for $Y$,
%To construct $X$ starting with $X_0=0$
 we need a different version of $X^{\langle 1\rangle}$, denoted by $\tilde X^{\langle 1\rangle}$, where given $\tilde X^{\langle 1\rangle}_n$, we set $\tilde X^{\langle 1\rangle}_{n+1}=\tilde X^{\langle 1\rangle}_n + v$ with probability $\omega^{\langle 1\rangle}_{\tilde X^{\langle 1\rangle}_n,\ell_n(\tilde X^{\langle 1\rangle}_n)}(v)$ if $\Pi_{d_1}v = 0$, and with probability $\omega^{\langle 1\rangle}_{\tilde X^{\langle 1\rangle}_n,\ell_n(\tilde X^{\langle 1\rangle}_n)}(v)/q(\Pi_{d_1}(v))$
if $\Pi_{d_1}v=Z_{N(\tilde X^{\langle 1\rangle}_{[0,n]})+1}$,
so that the increments of $\tilde X^{\langle 1\rangle}$ in the last $d_1$ coordinates match exactly with $(Z_n)_{n\geq 1}$. Let $(\tau^{\langle i\rangle}_j)_{j\geq 0}$ be defined for $X^{\langle i\rangle}$ as $(\tau_j)_{j\geq 0}$ is defined for $X$ in \refeq{tau}, and let $(\tilde \tau^{\langle 1\rangle}_j)_{j\geq 0}$ be defined similarly for $\tilde X^{\langle 1\rangle}$. Then we can construct
$(X_n)_{n\geq 0}$ by piecing together $(\tilde X^{\langle 1\rangle}_n)_{0\leq n\leq \tilde \tau^{\langle 1\rangle}_{T_1}}$ and $(X^{\langle i\rangle}_n)_{0\leq n\leq \tau^{\langle i\rangle}_{T_i-T_{i-1}}}$ for $i\geq 2$. More precisely, if we denote $\sigma_1=\tilde \tau^{\langle 1\rangle}_{T_1}$ and $\sigma_i = \sigma_{i-1} + \tau^{\langle i\rangle}_{T_i-T_{i-1}}$ for $i\geq 2$, then we set
\eq
\begin{aligned}
X_n  & = \ \tilde X^{\langle 1\rangle}_n  & \qquad \mbox{ for }\ & 0\leq n \leq \sigma_1, \\
X_n  & = \ X_{\sigma_1} + X^{\langle 2\rangle}_{n-\sigma_1} & \qquad \mbox{ for }\ & \sigma_1 \leq n \leq \sigma_2, \\
\vdots & & & \vdots \\
X_n  & = \ X_{\sigma_i} + X^{\langle i+1\rangle}_{n-\sigma_i} & \qquad \mbox{ for }\ & \sigma_i \leq n \leq \sigma_{i+1}, \\
\vdots & & & \vdots
\end{aligned}
\en
Conditional on $(Z_i)_{i\in\Z}$, $(X_n)_{n\geq 0}$ is thus represented as the concatenation of a sequence of random walks in independent cookie environments.  We leave it as an exercise to the reader to verify that (when averaged
over the law of $(Z_i)_{i\in\Z}$)
%, it is clear that
$X$ is distributed as a random walk in a random cookie environment $\omega$ with law $\P$.

To prove the law of large numbers, we write $X_{\tau_n} = \sum_{i=1}^n X_{\tau_i}-X_{\tau_{i-1}}$. Note that for each $n\in\Z$, there exists an $i\in\Z$ with $T_{i-1}<n\leq T_i$. We then define
$$
\Delta X^S_n := X^{\langle i\rangle}_{\tau^{\langle i\rangle}_{n-T_{i-1}}}-X^{\langle i\rangle}_{\tau^{\langle i\rangle}_{n-T_{i-1}-1}} \qquad \mbox{and} \qquad \Delta\tau^S_n := \tau^{\langle i\rangle}_{n-T_{i-1}} - \tau^{\langle i\rangle}_{n-T_{i-1}-1},
$$
which are the increments in space and time between successive jumps of $X^{\langle i\rangle}$ in the last $d_1$ coordinates.
Note that for each $i>T_1$, $\Delta X^S_i = X_{\tau_i}-X_{\tau_{i-1}}$ and $\Delta\tau^S_i=\tau_i-\tau_{i-1}$. The key to the proof of the law of large numbers is the observation that $(Z_i, \Delta X^S_i, \Delta\tau^S_i)_{i\in\Z}$ is an ergodic sequence with respect to the joint law of $(Z_i)_{i\in\Z}$ and $(X^{\langle i\rangle})_{i\in\Z}$. Assuming this, then by the ergodic theorem, almost surely,
$$
\lim_{n\to\infty} \frac{\sum_{i=1}^n \Delta X^S_i}{n} = \E[\Delta X^S_1]\in \R^d \qquad \mbox{ and } \qquad \lim_{n\to\infty} \frac{\sum_{i=1}^n \Delta\tau^S_i}{n} = \E[\Delta\tau^S_1] \in [1,\infty),
$$
where the ergodic theorem is applicable because $|\Delta X^S_1|_\infty \leq K\Delta\tau^S_1$, and $\Delta\tau^S_1$ is integrable because it is the time that a random
walk in the random environment $\omega^{\langle 1\rangle}$ has to wait before making a jump in the last $d_1$ coordinates, which by condition (b') is stochastically dominated by a geometric random variable with mean $\delta^{-1}$. Therefore, a.s.\ with respect to the law of $X$,
$$
\lim_{n\to\infty} \frac{X_{\tau_n}}{\tau_n} = \lim_{n\to\infty} \frac{\sum_{i=1}^n \Delta X^S_i}{\sum_{i=1}^n \Delta \tau^S_i} = \frac{\E[\Delta X^S_1]}{\E[\Delta\tau^S_1]} =: v \in \R^d.
$$
This implies \refeq{LLN} because $\sup_{\tau_{i-1}\leq n<\tau_i}|X_n-X_{\tau_{i-1}}|_\infty \leq K(\tau_i-\tau_{i-1})$ for each $i\in\N$, where $(\tau_i-\tau_{i-1})_{i\in\N}$
are dominated by independent geometric random variables with mean $\delta^{-1}$, and hence $\lim_{n\to\infty} (\tau_i-\tau_{i-1})/n\to 0$ almost surely.

It only remains to verify the ergodicity of $(Z_i, \Delta X^S_i, \Delta\tau^S_i)_{i\in\Z}$. Since $(Z_i)_{i\in\Z}$ determines the set of
cut times $\Di$, and conditional on $(Z_i)_{i\in\Z}$, $(\Delta X^S_n, \Delta \tau^S_n)_{T_{i-1}< n \leq T_n}$ is constructed independently for each $i\in\Z$
using only $(Z_n)_{T_{i-1}<n\leq T_i}$, by coupling, $(Z_i, \Delta X^S_i, \Delta\tau^S_i)_{i\in\Z}$ is shift invariant because
$(Z_i)_{i\in\Z}$ is shift invariant. The ergodicity of $(Z_i, \Delta X^S_i, \Delta\tau^S_i)_{i\in\Z}$ follows from the ergodicity of
the sequence $((Z_n)_{T_{i-1}<n\leq T_i}, (\Delta X^S_i)_{T_{i-1}<n\leq T_i}, (\Delta\tau^S_i)_{T_{i-1}<n\leq T_i})_{i\in\Z}$, the proof of which is
exactly the same as in the proof of \cite[Prop.~1.3]{BSZ03}.
\qed

\section{The lace expansion methodology}
\label{sec:lace-exp}
In this section we recall notation and results from \cite{HH07} and give a formula for the velocity $v$.

A nearest-neighbour
random walk path $\vec{x}_n$ is a sequence
$(x_i)_{i=0}^n$ for which $x_i=(x_i^{[1]},\dots,x_i^{[d]})\in \Z^d$
and $|x_{i+1}-x_i|=1$ (Euclidean distance) for each $i$.  If $\ewalk$ and $\walk$ are two such paths of length at least $j$ and $m$
respectively and such that $\eta_j=x_0$, then we define the concatenation
$\ewalk_j\circ \walk_m$ by
    \eq
    \lbeq{concat}
    (\ewalk_j\circ \walk_m)_{i}:=\left\{
    \begin{array}{lll}
    &\eta_i&\text{when }0\le i\leq j,\\
    &x_{i-j} &\text{when }j \leq i \leq m+j.
    \end{array}\right.
    \en

For a general nearest-neighbour path $\vec{x}_i$, we use the notation
$p^{\vec{x}_i}(x, y)$
for the conditional probability that the walk steps from $x$ (where $x=x_i$ is implicit in the notation) to
$y$, given the history of the path $\vec{x}_i=(x_0, \ldots, x_i)$.  In other words, for any finite path $\vec{x}_i$ of non-zero $P_{x_0}$ measure,
\eq
\lbeq{pdef}
p^{\vec{x}_i}(x_i, x_{i+1}):=P_{x_0}(X_{i+1}=x_{i+1}|\vec{X}_i=\vec{x}_i).
\en
Given $\ewalk_m$ such that $P_{\eta_0}(\vec{X}_m=\ewalk_m)>0$, we define a conditional probability measure $P^{\ewalk_m}$
on walks starting from $\eta_m$ by
    \eq
    \lbeq{cylinders}
    P^{\ewalk_m} (\vec{X}_n=\vec{x}_n)
    :=\prod_{i=0}^{n-1} p^{\ewalk_m\circ \vec{x}_i}(x_{i},x_{i+1})=P_{\eta_0}(\vec{X}_{m+n}=\ewalk_m\circ \vec{x}_n|\vec{X}_m=\ewalk_m).
    \en
Note that by definition, $P_{\eta_0}(\vec{X}_m=\walk_m)=P^{\eta_0}(\vec{X}_m=\walk_m)$.
%By the translation invariance of the model we have that $p^{\vec{\eta}_i}(x_i, x_{i+1})=p^{\vec{\eta}_i-\eta_0}(x_i-\eta_0, x_{i+1}-\eta_0)$.

Set $j_0=0$, and for $n\geq 1$ and fixed paths $\walkvec{n-1}{j_{n-1}+1}$ and $\walkvec{n}{j_{n}+1}$ let
    \eq
    \lbeq{Deltandef}
    \Delta_{\sss n}:= p^{\walkvec{n-1}{j_{n-1}+1}\circ
    \walkvec{n}{j_{n}}}(\walkcoor{n}{j_{n}},\walkcoor{n}{j_{n}+1})-p^{\walkvec{n}{j_{n}}}(\walkcoor{n}{j_{n}},\walkcoor{n}{j_{n}+1}),
    \en
which is a difference in the probabilities of stepping from $\walkcoor{n}{j_{n}}$ to $\walkcoor{n}{j_{n}+1}$ with two different histories,
$\walkvec{n-1}{j_{n-1}+1}\circ \walkvec{n}{j_{n}}$ and $\walkvec{n}{j_{n}}$, with the first history containing the second.

Define $\mc{A}_{m,\sN}:=\{(j_1, \dots, j_{\sN})\in \Z_+^{\sN}:\sum_{l=1}^Nj_l=m-N-1\}$, $\mc{A}_{\sN}:=\bigcupd_{m}\mc{A}_{m,\sN}=\Z_+^N$ and
      \eqalign
    \pi_m^{\smallsup{N}}(x,y):=&\sum_{\vec{j}\in \mc{A}_{m,N}}\sum_{\walkvec{0}{1}}\sum_{\walkvec{1}{j_{\sss 1}+1}}\dots\sum_{\walkvec{N}{j_{\sss N}+1} }1_{\{\eta^{(N)}_{j_{\sN}}=x, \eta^{(N)}_{j_{\sN}+1}=y\}}p^o(o,\walkcoor{0}{1})
    \prod_{n=1}^{N}\Delta_{\sss n}\prod_{i_{n}=0}^{j_{n}-1}p^{\walkvec{n-1}{j_{n-1}+1}\circ \walkvec{n}{i_{n}}}\left(\walkcoor{n}{i_{n}},\walkcoor{n}{i_{n}+1}\right),
    \lbeq{piNxydef}
    \enalign
where (here and throughout this paper), each $\sum_{\walkvec{i}{j_i+1}}$ is a sum over paths $(\walkcoor{i}{0},\dots,\walkcoor{i}{j_i+1})$ consisting of $j_i+1$ nearest-neighbour steps in $\Z^d$, where $\eta_0^{(0)}=o$ and otherwise $\walkcoor{i}{0}=\walkcoor{i-1}{j_{i-1}+1}$.  The summand is zero if the paths are not nearest-neighbour, so that we do not need to include this restriction in the summation notation. Note that $\pi_m^{\smallsup{N}}(x,y)=0$ for all $N\geq m$, and
by the indicator constraint in \refeq{piNxydef}, $\pi_m^{\smallsup{N}}(x,y)=0$ if $x$ is not a nearest-neighbour of $y$.  Furthermore, $\sum_y \pi_m^{\smallsup{N}}(x,y)=0$ since summing $\Delta_{\sss N}$ over $\walkcoor{N}{j_{N}+1}$ gives $1-1=0$ by \refeq{Deltandef}.

%The quantity \refeq{piNxydef} can also be expressed as
%\eq
%\pi_m^{\smallsup{N}}(x,y)=\sum_{\vec{j}\in \mc{A}_{m,N}}\sum_{\walkvec{0}{1}}\sum_{\walkvec{1}{j_{\sss %1}+1}}\dots\sum_{\walkvec{N}{j_{\sss N}+1} }I_{\{\omega^{(N)}_{j_{\sN}}=x, \omega^{(N)}_{j_{\sN}+1}=y\}}p^o(o,\walkcoor{0}{1})
%    \prod_{n=1}^{N}\Delta_{\sss n}P^{\walkvec{n-1}{j_{n-1}+1}}(S_{j_n}=\walkcoor{n}{j_n}).
%    \lbeq{piNxydef2}
%\en
%Note that this quantity depends on $\beta,\mu,d_*,d$, as do the following
Also define the following quantities
    \eq
    \pi_m(x,y):=\sum_{N=1}^{m-1} \pi_m^{\smallsup{N}}(x,y),  \quad
    \pi^{\smallsup{N}}(x,y):=\sum_{m=2}^\infty\pi_m^{\smallsup{N}}(x,y),
    \quad  \text{and} \quad \pi_m(y):=\sum_{N=1}^{m-1}\sum_x \pi_m^{\smallsup{N}}(x,y),
    \lbeq{piotherdef}
    \en
where (here and throughout this paper) $\sum_x$ denotes a sum over $x\in \Z^d$, where the summands are typically non-zero
only when $x$ is a nearest-neighbour of some $y\in\Z^d$ appearing in the summand.
%Note that the quantities $\pi_m^{\smallsup{N}}$ are all zero when $N+1>m$,
%and that
%All of the above quantities depend on $\beta$ and $\mu$.

The following result gives a formula for the velocity, provided the sum converges.
\begin{THM}[Proposition 3.1 of \cite{HH07}]
\label{thm:speedform}
Under the conditions of Theorem \ref{T:LLN}, the limiting velocity $v$ is given by
\eq
\lbeq{speedformula}
v=E_o[X_1] +\sum_{m=2}^\infty \sum_{y} y\pi_m(y),
\en
whenever this series converges.
\end{THM}
We are interested in properties of the first coordinate $v^{[1]}$ of the speed as a function of $\beta$.  Since $\sum_y \pi_m^{\smallsup{N}}(x,y)=0$ , \refeq{speedformula} can also be written in the more useful form
\eqalign
v^{}=E_o[X^{}_1] +\sum_{m=2}^\infty \sum_{N=1}^{m-1} \sum_{x,y} (y^{}-x^{})\pi_m^{\smallsup{N}}(x,y), \quad \text{ so that}\\
v^{[1]}=E_o[X^{[1]}_1] +\sum_{m=2}^\infty \sum_{N=1}^{m-1} \sum_{x,y} (y^{[1]}-x^{[1]})\pi_m^{\smallsup{N}}(x,y)\lbeq{speedformula2}.
\enalign

Under the assumptions of Theorem \ref{thm:monotone}, we will differentiate this expression with respect to $\beta$.  Note that
\eqalign
E_o[X^{[1]}_1]=&\mE[E_{\omega,o}[X_1^{[1]}]]=\mE[\omega_o(e_1)-\omega_o(-e_1)] \nonumber \\
=& \sum_{i=1}^2 (\nu_i(e_1)-\nu_i(-e_1))\mu_{d_*}(\nu_i) = \kappa\beta\rho+\kappa(\nu_1(e_1)-\nu_1(-e_1)), \nonumber
\enalign
where only the first term depends on $\beta$.  It follows immediately that
\eq
\frac{\del E_o[X^{[1]}_1]}{\del \beta}= \kappa\rho.\lbeq{mainderiv}
\en
If the derivative of the infinite series on the right hand side of \refeq{speedformula2} (with respect to $\beta$) is bounded in absolute value by $\kappa\rho$ then we will have shown that $v^{[1]}$ is increasing in $\beta$ (when $\rho>0$).  This is a strategy that has been used successfully in studying excited random walks \cite{HH09mono,H08compete}, where only in one coordinate direction does the step distribution not coincide with that of a simple symmetric random walk, and the probability of stepping in this one coordinate direction is a constant.

As one might infer from \refeq{piNxydef} and \refeq{Deltandef}, an analysis of the speed formula \refeq{speedformula}
and its derivative in $\beta$ reduces to an analysis of transition probabilities of the form \refeq{pdef}.

\subsection{The annealed transition probability}
\label{sec-transprob}
In this section we consider properties of the annealed transition probability, defined for any path $\vec{\eta}_n$ of positive $P_{\eta_0}$-measure by
\eqalign
\lbeq{anntransprob}
p^{\vec{\eta}_n}(\eta_n,\eta_{n+1}):=&P_{\eta_0}(X_{n+1}=\eta_{n+1}|\vec{X}_n=\vec{\eta}_n)=\frac{P_{\eta_0}(\vec{X}_{n+1}=\vec{\eta}_{n+1})}{P_{\eta_0}(\vec{X}_n=\vec{\eta}_n)}\nn\\
=&\frac{\E[P_{\omega,\eta_0}(\vec{X}_{n+1}=\vec{\eta}_{n+1})]}{\E[P_{\omega,\eta_0}(\vec{X}_{n}=\vec{\eta}_{n})]}=
%\frac{\int P_{\omega,\eta_0}(X_{n}=\eta_{n})\omega_{\eta_n}(\eta_{n+1}-\eta_n) d\mu}{\int P_{\omega,\eta_0}(X_{n}=\eta_{n})d\mu}\nn
%=&
\frac{\E\left[\omega_{\eta_n}(\eta_{n+1}-\eta_{n})\prod_{i=0}^{n-1} \omega_{\eta_{i}}(\eta_{i+1}-\eta_{i})\right]}{\E[ \prod_{i=0}^{n-1} \omega_{\eta_{i}}(\eta_{i+1}-\eta_{i})]}.
\enalign
Under $\mP$, $\omega_x$ and $\omega_y$ are independent if $x\ne y$, whence \refeq{anntransprob} is equal to
\eqalign
\lbeq{anntransprob_2}
%\frac{
\frac{\mE[\omega_{\eta_n}(\eta_{n+1}-\eta_{n})B(\vec{\eta}_{n})]}{
%}{
\mE[B(\vec{\eta}_{n})]},
\enalign
where
\[B(\vec{\eta}_{n}):=\prodtwo{0\le j\le n-1:}{\eta_j=\eta_n} \omega_{\eta_{n}}(\eta_{j+1}-\eta_{j}).\]
Therefore
\eq
\lbeq{Delta2}
\Delta_{\sss n} =\frac{\mE[\omega_{\walkcoor{n}{j_n}}(\walkcoor{n}{j_n+1}-\walkcoor{n}{j_n})B(\walkvec{n-1}{j_{n-1}+1}\circ \walkvec{n}{j_n})]}{\mE[B(\walkvec{n-1}{j_{n-1}+1}\circ \walkvec{n}{j_n})]}
-\frac{\mE[\omega_{\walkcoor{n}{j_n}}(\walkcoor{n}{j_n+1}-\walkcoor{n}{j_n})B(\walkvec{n}{j_n})]}{\mE[B(\walkvec{n}{j_n})]}.
\en
It follows immediately that $\Delta_n\neq 0$ only if $B(\walkvec{n}{j_n})\ne B(\walkvec{n-1}{j_{n-1}+1}\circ \walkvec{n}{j_n})$, i.e.~only if $\walkcoor{n}{j_n}\in \walkvec{n-1}{j_{n-1}}$.
% we thus have
%\eqalign
%\lbeq{Delta_n_bound}
%|\Delta_{\sss n}|\le &1_{\{\eta^{(n)}_{j_{\sss n}}\in \walkvec{n-1}{j_{n-1}}\}}\begin{cases}
%(1-\delta)q(\walkcoor{n}{j_n+1}-\walkcoor{n}{j_n}), &\text{ if }\walkcoor{n}{j_n+1}-\walkcoor{n}{j_n}\in V_d\backslash V_{d_0}, \\
%1-\delta, &\text{ if }\walkcoor{n}{j_n+1}-\walkcoor{n}{j_n}\in V_{d_0},
%\end{cases}
%\enalign
Recall that $V_k:=\{\pm e_i\}_{1\leq i\leq k}$.
By Assumptions {\bf (A0)-(A4)}, for $\P$-a.s.\ all $\omega_o$, $\delta q(v)\leq \omega_o(v)\leq q(v)$ for all $v\in V_d\backslash V_{d_0}$.
Therefore $\sum_{v\in V_{d_0}}\omega_o(v)\leq 1-\delta$, and
\eq\lbeq{deltanbd}
\sum_{\eta^{(n)}_{j_{\sss n}+1}\in\Z^d} |\Delta_{\sss n}| \leq 2(1-\delta)1_{\{\eta^{(n)}_{j_{\sss n}}\in \walkvec{n-1}{j_{n-1}}\}},
\en
since the two terms in \refeq{Delta2} represent two probability kernels on $\Z^d$ which both dominate $\delta q(\cdot)$.

We also need to examine the derivatives of the annealed transition probabilities with respect to $\beta$.
For a directed edge $b$, let $\ell(\vec{\eta}_n,b)=\sum_{i=1}^n 1_{\{(\eta_{i-1},\eta_i)=b\}}$
denote the edge local time of $\vec{\eta}$ at $b$ up to time $n$, and for any $V\subset V_d$ let $\ell(\vec{\eta}_n,V)=\sum_{b \in V}\ell(\vec{\eta}_n,(\eta_n,\eta_n+b))$.  Then for each $\vec{\eta}_n$, almost surely, at most one of the following can be greater than 0:
\eqalign \lbeq{edgeloctime}
\ell(\vec{\eta}_n,\mc{S}_1), \quad \ell(\vec{\eta}_n,\mc{S}_2), \quad \ell(\vec{\eta}_n,V_{d_*}\setminus(\mc{S}_1\cup \mc{S}_2)),
\enalign
where we recall from {\bf (A2)} that $\mc{S}_i={\rm supp}(\nu_i)\subset V_{d_*}$, $i=1,2$.
\begin{itemize}
\item If $u_{n}:=\eta_{n+1}-\eta_n\in \mc{S}_1\cup \mc{S}_2$, then
\eqalign
p^{\vec{\eta}_n}(\eta_n,\eta_{n+1})=&
\sum_{i=1}^2 \nu_i(u_n)1_{\{\ell(\vec{\eta}_n,\mc{S}_i)>0\}}+1_{\{\ell(\vec{\eta}_n,V_{d_*})=0\}}\sum_{i=1}^2\nu_i(u_n)\mu_{d_*}(\nu_i) \nn\\
=&\sum_{i=1}^2 \nu_i(u_n)\left[1_{\{\ell(\vec{\eta}_n,\mc{S}_i)>0\}}+1_{\{\ell(\vec{\eta}_n,V_{d_*})=0\}}\mu_{d_*}(\nu_i)\right],\nn
\enalign
from which we deduce
\eqalign
\lbeq{trans_deriv}
\frac{\del p^{\vec{\eta}_n}(\eta_n,\eta_{n+1})}{\del \beta}=&
\kappa[\nu_2(u_n)-\nu_1(u_n)]1_{\{\ell(\vec{\eta}_n,V_{d_*})=0\}}.
\enalign
\item If $u_{n}:=\eta_{n+1}-\eta_n\notin \mc{S}_1\cup \mc{S}_2$, then it is easily verified by direct calculations that
\eqalign
\frac{\del p^{\vec{\eta}_n}(\eta_n,\eta_{n+1})}{\del \beta}=0. \nn
\enalign
\end{itemize}
Therefore
\eqalign
\lbeq{Delta_n_deriv}
\frac{\del \Delta_{\sss n}}{\del \beta}=&
\kappa[\nu_2(\walkcoor{n}{j_{n}+1}-\walkcoor{n}{j_{n}})-\nu_1(\walkcoor{n}{j_{n}+1}-\walkcoor{n}{j_{n}})][1_{\{\ell(\walkvec{n-1}{j_{n-1}+1}\circ \walkvec{n}{j_{n}},V_{d_*})=0\}}-1_{\{\ell(\walkvec{n}{j_{n}},V_{d_*})=0\}}],
\enalign
so that
\eqalign
\nn
\left|\frac{\del \Delta_{\sss n}}{\del \beta}\right|\le &
\kappa|\nu_2(\walkcoor{n}{j_{n}+1}-\walkcoor{n}{j_{n}})-\nu_1(\walkcoor{n}{j_{n}+1}-\walkcoor{n}{j_{n}})|1_{\{\walkcoor{n}{j_{n}}\in \walkvec{n-1}{j_{n-1}}\}},
\enalign
which together with the fact $\nu_1, \nu_2\in \Mi_\gamma(\Z^{d_*})$ and $\gamma\leq 1-\delta$ implies that
\eqalign
\lbeq{Delta_n_deriv_bound}
\sum_{\walkcoor{n}{j_{n}+1}\in\Z^d} \left|\frac{\del \Delta_{\sss n}}{\del \beta}\right|\le & 2\kappa(1-\delta)1_{\{\walkcoor{n}{j_{n}}\in \walkvec{n-1}{j_{n-1}}\}}.
\enalign
Observe that
\eqalign
&\sum_y(y^{[1]}-x^{[1]})\big(p^{\vec{\eta}_n}(x,y)-p^{\vec{x}_m\circ\vec{\eta}_n}(x,y)\big)\nn\\
\qquad =&1_{\{x\in \vec{x}_{m-1}\}}[p^{\vec{\eta}_n}(x,x+e_1)-p^{\vec{x}_m\circ\vec{\eta}_n}(x,x+e_1)-p^{\vec{\eta}_n}(x,x-e_1)+p^{\vec{x}_m\circ\vec{\eta}_n}(x,x-e_1)].\lbeq{pdrift}
\enalign
\begin{LEM}
\label{lem:driftdelta}
For all $x\in \Z^d$ and nearest-neighbour paths $\vec{x}_m$ and $\vec{\eta}_n$ such that $\eta_0=x_m$,
\eqalign
|\sum_y(y^{[1]}-x^{[1]})\big(p^{\vec{\eta}_n}(x,y)-p^{\vec{x}_m\circ\vec{\eta}_n}(x,y)\big)|\le &\ \rho 1_{\{x\in \vec{x}_{m-1}\}}, \lbeq{Delta_N_bound}\\
\big|\frac{\del}{\del \beta}\sum_y(y^{[1]}-x^{[1]})\big(p^{\vec{\eta}_n}(x,y)-p^{\vec{x}_m\circ\vec{\eta}_n}(x,y)\big)\big|\le &\ \kappa \rho  1_{\{x\in \vec{x}_{m-1}\}}. \lbeq{Delta_N_deriv_bound}
\enalign
\end{LEM}
{\bf Proof.}
The term in brackets on the right hand side of \refeq{pdrift} is equal to
\eqalign
&\sum_{i=1}^2 \nu_i(e_1)\left[1_{\{\ell(\vec{\eta}_n,\mc{S}_i)>0\}}+1_{\{\ell(\vec{\eta}_n,V_{d_*})=0\}}\mu_{d_*}(\nu_i)\right]\nn\\
&\quad\quad -\sum_{i=1}^2 \nu_i(e_1)\left[1_{\{\ell(\vec{x}_m\circ\vec{\eta}_n,\mc{S}_i)>0\}}+1_{\{\ell(\vec{x}_m\circ\vec{\eta}_n,V_{d_*})=0\}}\mu_{d_*}(\nu_i)\right]\nn\\
&\quad\quad\qquad -\sum_{i=1}^2 \nu_i(-e_1)\left[1_{\{\ell(\vec{\eta}_n,\mc{S}_i)>0\}}+1_{\{\ell(\vec{\eta}_n,V_{d_*})=0\}}\mu_{d_*}(\nu_i)\right]\nn\\
&\quad\quad\qquad\qquad +\sum_{i=1}^2 \nu_i(-e_1)\left[1_{\{\ell(\vec{x}_m\circ\vec{\eta}_n,\mc{S}_i)>0\}}+1_{\{\ell(\vec{x}_m\circ\vec{\eta}_n,V_{d_*})=0\}}\mu_{d_*}(\nu_i)\right]\nn\\
%=&\sum_{i=1}^2 [\nu_i(e_1)-\nu_i(-e_1)]\left[1_{\{\ell(\vec{\eta}_n,\mc{S}_i)>0\}}+1_{\{\ell(\vec{\eta}_n,V_{d_*})=0\}}\mu_{d_*}(\nu_i)\right]\\
%&-\sum_{i=1}^2 [\nu_i(e_1)-\nu_i(-e_1)]\left[1_{\{\ell(\vec{x}_m\circ\vec{\eta}_n,\mc{S}_i)>0\}}+1_{\{\ell(\vec{x}_m\circ\vec{\eta}_n,V_{d_*})=0\}}\mu_{d_*}(\nu_i)\right]\\
=&\sum_{i=1}^2 [\nu_i(e_1)-\nu_i(-e_1)]\Big[1_{\{\ell(\vec{\eta}_n,\mc{S}_i)>0\}}+1_{\{\ell(\vec{\eta}_n,V_{d_*})=0\}}\mu_{d_*}(\nu_i)\nn\\
&\qquad \qquad\qquad\qquad\qquad\qquad -1_{\{\ell(\vec{x}_m\circ\vec{\eta}_n,\mc{S}_i)>0\}}-1_{\{\ell(\vec{x}_m\circ\vec{\eta}_n,V_{d_*})=0\}}\mu_{d_*}(\nu_i)\Big].\nn
\enalign
If the first indicator function is non-zero for some $i$ then so is the third (for the same $i$), while all other indicators are zero.  Therefore we can rewrite the above as
$$
\sum_{i=1}^2 [\nu_i(e_1)-\nu_i(-e_1)]1_{\{\ell(\vec{\eta}_n,\mc{S}_1\cup \mc{S}_2)=0\}}\Big[1_{\{\ell(\vec{\eta}_n,V_{d_*})=0\}}\mu_{d_*}(\nu_i)
-1_{\{\ell(\vec{x}_m\circ\vec{\eta}_n,\mc{S}_i)>0\}}-1_{\{\ell(\vec{x}_m\circ\vec{\eta}_n,V_{d_*})=0\}}\mu_{d_*}(\nu_i)\Big].
$$
If the final indicator function here is 1 then so is the second, while the third is zero.  Thus the quantity above is zero unless the last indicator is zero, i.e.~\refeq{pdrift} is equal to
\eqalign
&1_{\{x \in \vec{x}_{m-1}\}}1_{\{\ell(\vec{\eta}_n,\mc{S}_1\cup \mc{S}_2)=0\}}1_{\{\ell(\vec{x}_m\circ\vec{\eta}_n,V_{d_*})>0\}}\nn\\
&\qquad \qquad \qquad \qquad \qquad \times \sum_{i=1}^2 [\nu_i(e_1)-\nu_i(-e_1)]\Big[1_{\{\ell(\vec{\eta}_n,V_{d_*})=0\}}\mu_{d_*}(\nu_i)-1_{\{\ell(\vec{x}_m\circ\vec{\eta}_n,\mc{S}_i)>0\}}\Big]\lbeq{pdrift2}.
\enalign
\begin{itemize}
\item[Case 1:] If $\nu_1(e_1)=\nu_1(-e_1)=0$, then $\rho=\nu_2(e_1)-\nu_2(-e_1)$, and \refeq{pdrift2} becomes
\eqalign
&1_{\{x \in \vec{x}_{m-1}\}}1_{\{\ell(\vec{\eta}_n,\mc{S}_1\cup \mc{S}_2)=0\}}1_{\{\ell(\vec{x}_m\circ\vec{\eta}_n,V_{d_*})>0\}}\;\rho\;\Big[1_{\{\ell(\vec{\eta}_n,V_{d_*})=0\}}\kappa\beta-1_{\{\ell(\vec{x}_m\circ\vec{\eta}_n,\mc{S}_2)>0\}}\Big],\nn
\enalign
where the term in brackets is the difference of two terms between 0 and 1 and hence is bounded in absolute value by 1.  The derivative with respect to $\beta$ is
\eq
1_{\{x \in \vec{x}_{m-1}\}}1_{\{\ell(\vec{\eta}_n,\mc{S}_1\cup \mc{S}_2)=0\}}1_{\{\ell(\vec{x}_m\circ\vec{\eta}_n,V_{d_*})>0\}}1_{\{\ell(\vec{\eta}_n,V_{d_*})=0\}}\;\kappa\; \rho \nn
\en
which is bounded in absolute value by $1_{\{x \in \vec{x}_{m-1}\}}\kappa\rho$.
\item[Case 2:] If $\nu_1(e_1)=0$ and $\nu_1(-e_1)>0$, then $\nu_2(-e_1)=0$ and $\rho=\nu_2(e_1)+\nu_1(-e_1)$, while
\refeq{pdrift2} becomes
\eqalign
&1_{\{x \in \vec{x}_{m-1}\}}1_{\{\ell(\vec{\eta}_n,\mc{S}_1\cup \mc{S}_2)=0\}}1_{\{\ell(\vec{x}_m\circ\vec{\eta}_n,V_{d_*})>0\}}\times\nn\\
& \Bigg[\nu_2(e_1)\Big[1_{\{\ell(\vec{\eta}_n,V_{d_*})=0\}}\kappa\beta-1_{\{\ell(\vec{x}_m\circ\vec{\eta}_n,\mc{S}_2)>0\}}\Big]-\nu_1(-e_1)\Big[1_{\{\ell(\vec{\eta}_n,V_{d_*})=0\}}\kappa(1-\beta)-1_{\{\ell(\vec{x}_m\circ\vec{\eta}_n,\mc{S}_1)>0\}}\Big]\Bigg],\nn
\enalign
where the term inside the largest brackets is bounded in absolute value by $\rho$.  The derivative with respect to $\beta$ is
\eq
1_{\{x \in \vec{x}_{m-1}\}}1_{\{\ell(\vec{\eta}_n,\mc{S}_1\cup \mc{S}_2)=0\}}1_{\{\ell(\vec{x}_m\circ\vec{\eta}_n,V_{d_*})>0\}}1_{\{\ell(\vec{\eta}_n,V_{d_*})=0\}}\;\kappa\;[\nu_2(e_1)+\nu_1(-e_1)],\nn
\en
which is bounded in absolute value by $1_{\{x \in \vec{x}_{m-1}\}}\kappa \rho$.
\item[Case 3:] If $\nu_1(e_1)>0$ and $\nu_1(-e_1)=0$, then $\nu_2(e_1)=0$ and $\rho=-(\nu_2(-e_1)+\nu_1(e_1))<0$, contradicting our assumption $\rho>0$
in {\bf (A2)}. So this case can be ruled out.
%while \refeq{pdrift2} equals
%\eqalign
%&1_{\{x \in \vec{x}_{m-1}\}}1_{\{\ell(\vec{\eta}_n,\mc{S}_1\cup \mc{S}_2)=0\}}1_{\{\ell(\vec{x}_m\circ\vec{\eta}_n,V_{d_*})>0\}}\nn\\
%&\times %\Bigg[\nu_1(e_1)\Big[1_{\{\ell(\vec{\eta}_n,V_{d_*})=0\}}\kappa(1-\beta)-1_{\{\ell(\vec{x}_m\circ\vec{\eta}_n,\mc{S}_1)>0\}}\Big]-\nu_2(-e_1)\Big[1_{\{\ell(\vec{\eta}_n,V_{d_*})=0\}}\kappa\beta-1_{\{\ell(\vec{x}_m\circ\vec{\eta}_n,\mc{S}_2)>0\}}\Big]\Bigg],\nn
%\enalign
%where the term inside the largest brackets is bounded in absolute value by $|\rho|=\nu_2(-e_1)+\nu_1(e_1)$.  The derivative with respect to $\beta$ is
%\eq
%-1_{\{x \in \vec{x}_{m-1}\}}1_{\{\ell(\vec{\eta}_n,\mc{S}_1\cup %\mc{S}_2)=0\}}1_{\{\ell(\vec{x}_m\circ\vec{\eta}_n,V_{d_*})>0\}}1_{\{\ell(\vec{\eta}_n,V_{d_*})=0\}}\;\kappa\;[\nu_1(e_1)+\nu_2(-e_1)],\nn
%\en
%which is bounded in absolute value by $1_{\{x \in \vec{x}_{m-1}\}}\kappa|\rho|$.
\item [Case 4:] If $\nu_1(e_1)>0$ and $\nu_1(-e_1)>0$, then $\nu_2(e_1)=\nu_2(-e_1)=0$ and $\rho=\nu_1(-e_1)-\nu_1(e_1)$, while \refeq{pdrift2} equals
\eqalign
&-1_{\{x \in \vec{x}_{m-1}\}}1_{\{\ell(\vec{\eta}_n,\mc{S}_1\cup \mc{S}_2)=0\}}1_{\{\ell(\vec{x}_m\circ\vec{\eta}_n,V_{d_*})>0\}}\;\rho\;\Big[1_{\{\ell(\vec{\eta}_n,V_{d_*})=0\}}\kappa(1-\beta)-1_{\{\ell(\vec{x}_m\circ\vec{\eta}_n,\mc{S}_1)>0\}}\Big],\nn
\enalign
with the term inside the bracket bounded in absolute value by $1$. The derivative with respect to $\beta$ is
\eq
1_{\{x \in \vec{x}_{m-1}\}}1_{\{\ell(\vec{\eta}_n,\mc{S}_1\cup \mc{S}_2)=0\}}1_{\{\ell(\vec{x}_m\circ\vec{\eta}_n,V_{d_*})>0\}}1_{\{\ell(\vec{\eta}_n,V_{d_*})=0\}}\,\kappa\,\rho,\nn
\en
which is bounded in absolute value by $1_{\{x \in \vec{x}_{m-1}\}}\kappa\rho$.
\end{itemize}
This completes the proof of the lemma.\qed

\subsection{Convergence of the speed formula}
\label{sec:formula_converges}
In this section we prove a bound on $\sum_{x, y}\sum_m|\pi^{\smallsup{N}}_m(x,y)|$, which can then be used to show that the sum in the speed formula \refeq{vformula}  (or more precisely, the first line of \refeq{speedformula2})
converges, and hence \refeq{vformula} holds. Similar methods have been used in \cite{HH09mono} and \cite{H08compete} to bound similar quantities.
The present context is more demanding since the types of environment being considered are more complicated (and require a somewhat more detailed analysis).  In addition, the probability of stepping in various coordinate directions is allowed to be random, so that the number of steps taken in the coordinate directions $d_0+1,\dots,d$ is not binomially distributed (as in \cite{HH09mono} and \cite{H08compete}), but rather is stochastically dominated by a binomial distribution.
%For readers with no interest in explicit estimates, the important point is that in Proposition \ref{prp:pibound} the bounds are finite and summable when $d_1\ge 5$ and $\delta$ is sufficiently small, and converge to zero as $\delta\searrow 0$.

%For $i\ge 0$, define {\ch not needed any more!!!}
%    \eqalign
%\mc{E}_i=\mc{E}_i(q,\delta)=&
%   \sup_{v \in \Z^{d_1}}\big(\delta^{-(i+1)}G^{*(i+1)}_{q}(v)-1_{\{v=o\}}\big).\lbeq{Eidef}
%    \enalign

We will need the following extension of \cite[Lemma 3.1]{HH09mono} for our bounds.

\begin{LEM}
\label{lem:togreens}
Let $\vec{X}$ be a random walk in $\Z^d=\Z^{d_0+d_1}$ in a random environment $\omega$ satisfying assumptions {\bf (A0)--(A4)}. For any $\vec\eta_m$ with
$P_{\eta_0}(\vec X_m=\vec\eta_m)>0$, $i \in \Z_+$ and $u\in \Z^d$, if we denote $P^{\vec \eta_m}(\cdot):=P_{\eta_0}(\cdot |\,(X_{-m+k})_{0\leq k\leq m}=\vec\eta_m)$, then
%\begin{itemize}
%\item if $P^{\vec \eta_n}(\cdot):=P_{\eta_0}(\cdot |\vec X_n=\vec\eta_n\big)$, then there exists $\delta\in (0,1)$ such that $P^{\vec \eta_n}\big(X_{n+1}-X_n\notin V_{d_0})\ge \delta$ for all $n$ and $\vec{\eta}_n$ and
%\item the sequence of steps in the $d_1$ coordinate directions constitute simple random walk with step distribution $q(\cdot)$,
% such that $G_q^{*2}<\infty$
%\end{itemize}
\eqalign
    \sum_{j=0}^{\infty}\frac{(j+i)!}{j!}P^{\vec{\eta}_m}(X_{j}=u) & \ \le\ i!\delta^{-(i+1)}G^{*(i+1)}_q,\lbeq{Gbound1}
%{\bf \text{(no longer needed)}}\quad     &\sum_{j=1}^{\infty}\frac{(j+i)!}{j!}P^{\vec{\eta}_m}(X_{j}=u)
%    \le i!\mc{E}_i,\lbeq{Gbound2}\\
%    \sum_{j=1}^{\infty}\sum_{l=0}^{j-1}\frac{(j+i)!}{j!}P^{\vec{\eta}_m}(X_{j}=u|X_{l+1}-X_l=w)
%    & \ \le\  (i+1)!\delta^{-(i+2)}G_q^{*(i+2)},\lbeq{Gbound3}
    \enalign
where $G^{*k}_q$ is defined in (\ref{Gcond}).
%{\bf \text{(no longer needed)}}\quad       \text{ for any }w\in V_{d_0}, \quad &\sum_{j=2}^{\infty}\sum_{l=0}^{j-1}\frac{(j+i)!}{j!}P^{\vec{\eta}_m}(X_{j}=u|X_{l+1}-X_l=w)
%    \le (i+1)!\mc{E}_{i+1}\lbeq{Gbound4}.
\end{LEM}
\noindent
{\bf Proof.}
Let $\mc{N}_j$ be the number of steps the walk $\vec X_j:=(X_{k})_{0\leq {k}\leq j}$ takes in the last $d_1$ coordinate directions, given history $(X_{-m+{k}})_{0\leq { k}\leq m}=\vec{\eta}_m$. Let $\tau_n:=\inf\{j\geq 0: \mc{N}_j=n\}$. Let $(Y_n)_{n\in\Z_+}$ be the random walk on $\Z^{d_1}$ coupled with $X$ such that $Y_0=\Pi_{d_1}(X_0)$,
the projection of $X_0\in\Z^{d_0+d_1}$ to its last $d_1$ coordinates, and $Y_n-Y_{n-1}=\Pi_{d_1}(X_{\tau_n}-X_{\tau_{n-1}})$ for all $n\in\N$. By assumption {\bf (A0)},
it is clear that $Y$ is distributed as a random walk on $\Z^{d_1}$ with transition kernel $q$. We will denote the law of $Y$ separately by $\mc{P}_q$. We then have
\eqalign\lbeq{green1}
P^{\vec{\eta}_m}(X_j=u)=\sum_{l=0}^jP^{\vec{\eta}_m}(X_j=u,\mc{N}_j=l)\le& \sum_{l=0}^j P^{\vec{\eta}_m}(\Pi_{d_1}(X_j)=\Pi_{d_1}(u),\mc{N}_j=l)\nn\\
%=&\sum_{l=0}^n P^{\vec{\eta}_m}(Y_l=\Pi_{d_1}(u),\mc{N}_n=l)\nn\\
=&\sum_{l=0}^j P^{\vec{\eta}_m}(\mc{N}_j=l | Y_l=\Pi_{d_1}(u)) \mc{P}_q(Y_l=\Pi_{d_1}(u)).
\enalign
Combined with Lemma \ref{lem:qpart2} below, we obtain
\begin{eqnarray*}
\sum_{j=0}^{\infty}\frac{(j+i)!}{j!}P^{\vec{\eta}_m}(X_{j}=u) &\leq& \sum_{j=0}^\infty \frac{(j+i)!}{j!} \sum_{l=0}^j P^{\vec{\eta}_m}(\mc{N}_j=l | Y_l=\Pi_{d_1}(u)) \mc{P}_q(Y_l=\Pi_{d_1}(u)) \\
&=& \sum_{l=0}^\infty \mc{P}_q(Y_l=\Pi_{d_1}(u)) \sum_{j=0}^{\infty}\frac{(j+i)!}{j!} P^{\vec{\eta}_m}(\mc{N}_j=l | Y_l=\Pi_{d_1}(u)) \\
&\leq& \sum_{l=0}^\infty \delta^{-i}\frac{(l+i)!}{l!} \mc{P}_q(Y_l=\Pi_{d_1}(u)).
\end{eqnarray*}
The inequality \refeq{Gbound1} then follows from the fact that (see e.g.~\cite[(3.2)]{HH09mono})
\eq\lbeq{Gform}
G^{*(i+1)}_q(v)=\sum_{l=0}^{\infty}\frac{(l+i)!}{i!\,l!}\mc{P}_q(Y_l=v).
\en
\qed

\begin{LEM}
\label{lem:qpart2} Let $\vec X$, $\vec{\eta}_m$, $u$, $\mc{N}_j$ and $\vec Y$ be defined as in Lemma \ref{lem:togreens} and its proof. Then
\eq
\sum_{j=0}^{\infty}\frac{(j+i)!}{j!}P^{\vec{\eta}_m}(\mc{N}_j=l|Y_l=\Pi_{d_1}(u))\le \delta^{-i}\frac{(l+i)!}{l!}.
\en
\end{LEM}
\noindent
{\bf Proof.} First we claim that for all $\vec x_k$ and $\vec y_l$ which are compatible (by the coupling of $\vec X$ and $\vec Y$), we have
\eq
P^{\vec{\eta}_m}(\Pi_{d_1}(X_{k+1}-X_k)\neq 0|\vec{Y}_l=\vec{y}_l,\vec{X}_k=\vec{x}_k)\ge \delta.\lbeq{crucial}
\en
Note that
\eqalign \nn
\sum_{j=0}^{\infty}\frac{(j+i)!}{j!}P^{\vec{\eta}_m}(\mc{N}_j=l|Y_l=\Pi_{d_1}(u))=&E^{\vec{\eta}_m}\big[\sum_{j=0}^{\infty}\frac{(j+i)!}{j!}1_{\{\mc{N}_j=l\}}\big|Y_l=\Pi_{d_1}(u)\big]\\
=&E^{\vec{\eta}_m}\big[\sum_{j=\tau_l}^{\tau_{l+1}-1}\frac{(j+i)!}{j!}\big|Y_l=\Pi_{d_1}(u)\big],\lbeq{to_couple}
\enalign
where $\tau_i=\tau_i(\mc{N})$ is the first hitting time of level $i$ by $\mc{N}_j$.  By \refeq{crucial}, regardless of $\vec{Y}_l,\vec{X}_k$ and $\vec{\eta}_m$, the $(k+1)$st step has probability at least $\delta$ of having a non-zero $\Pi_{d_1}$ projection. It follows that under any conditional measure depending only on $\vec{Y}_l$ and $\vec{\eta}_m$, we can couple $\mc{N}_j$ with a random walk $\mc{M}_j$ on $\Z_+$ taking i.i.d.~steps +1 or 0 with probabilities $\delta$ and $1-\delta$ respectively, such that $\tau_{i+1}(\mc{N})-\tau_i(\mc{N})\le \tau_{i+1}(\mc{M})-\tau_i(\mc{M})$ for all $i$, almost surely.  This also implies that $\mc{M}_j\le \mc{N}_j$ and $\tau_i(\mc{N})\le \tau_i(\mc{M})$ a.s.

Note that $\mc{M}_j$ is binomial with parameters $(j,\delta)$.  Therefore \refeq{to_couple} is bounded by
\eqalign\lbeq{comb1}
E\Big[\sum_{j=\tau_l(\mc{M})}^{\tau_{l+1}(\mc{M})-1}\frac{(j+i)!}{j!}\Big]=\sum_{j=l}^{\infty}\frac{(j+i)!}{j!}P(\mc{M}_j=l)=\sum_{j=l}^{\infty}\frac{(j+i)!}{j!}{j\choose l}\delta^{l}(1-\delta)^{j-l}=\delta^{-i}\frac{(l+i)!}{l!},
\enalign
exactly as in \cite{HH09mono,H08compete}.  It therefore remains to prove \refeq{crucial}.

To prove \refeq{crucial}, recall \refeq{tau} and let $J_{k}=\max\{i\ge 0:\tau_i\le k\}$. Then $\vec{Y}_{J_{k}}$ is determined by
$\vec{X}_{k}$ by the coupling of $\vec X$ and $\vec Y$, and by assumptions (b') and {\bf (A0)}, $(Y_t-Y_{J_{k}})_{t\geq J_{k}}$ is independent of $\vec X_k$ and the
event $\{\Pi_{d_1}(X_{k+1}-X_{k})\ne 0\}$. Therefore
\eqalign
P^{\vec{\eta}_m}(\Pi_{d_1}(X_{k+1}-X_k)\neq 0|\vec{Y}_l=\vec{y}_l,\vec{X}_k=\vec{x}_k)=&P^{\vec{\eta}_m}(\Pi_{d_1}(X_{k+1}-X_{k})\neq 0|\vec{X}_{k}=\vec{x}_k)\nn\\
=& E^{\vec{\eta}_m} \big[\sum_{u: \Pi_{d_1}(u)\neq 0}\omega_{x_k}(u)\big|\vec{X}_{k}=\vec{x}_k\big]\ge \delta,\nn
\enalign
%Note that conditional on $X_{k+m}$ and $\omega_{X_{k+m}}$, the law of $Y$ remains the same. Therefore
%, the expectation above equals
%\eq
%E\big[\sum_{u: \Pi_{d_1}(u)\neq 0}\omega_{X_{k+m}}(u)\big|\vec{X}_{k+m}=\vec{\eta}_m\circ\vec{x}_k\big] \ge \delta,
%\en
as required, where we used assumption (b') once more in the inequality.
\qed

\vspace{.5cm}

Let
\eq
\lbeq{alphadef}
\epsilon_\delta:=2(1-\delta)\qquad \mbox{and} \qquad \alpha=\epsilon_{\delta}\delta^{-2}G^{*2}_q.
\en
The following proposition, together with Proposition 3.1 of \cite{HH07}, shows that the series in the speed formula \refeq{vformula} converges when $G_q^{*2}<\infty$ and $\alpha<1$.  When $G_q^{*2}<\infty$ as is assumed in {\bf (A3)}, $\alpha<1$ can be achieved if $\delta<1$ is sufficiently close to $1$.
\begin{PRP}
\label{prp:pibound}
For RWpRE as in Theorem \ref{thm:monotone} and for each $N\in\N$, we have
\eqalign
\sum_{x,y\in\Z^d}\sum_{m=2}^\infty |\pi_{m}^{\smallsup{N}}(x,y)|\le &
\epsilon_{\delta}\delta^{-1}G_{q} \alpha^{N-1}.
 \lbeq{pibound}
\enalign
\end{PRP}
%It then follows from Theorem 1.4 of \cite{BSZ03}, Proposition 3.1 of \cite{HH07} and the fact that $G_5^{*2}<\infty$, that \refeq{speedformula2} holds provided $d_1\ge 5$ and $\alpha
 %<1$.
\noindent
{\bf Proof.}
It follows from \refeq{piNxydef} that $\sum_{x,y\in\Z^d}\sum_{m=2}^\infty|\pi_m^{\smallsup{N}}(x,y)|$ is bounded by
%and \refeq{Deltabound} that
\eqalign
&\sum_{\walkcoor{0}{1}}p^{o}(o,\walkcoor{0}{1})
\sum_{j_{\sss 1}=1}^{\infty}\sum_{\walkvec{1}{j_{\sss 1}+1}}
|\Delta_{\sss 1}|\prod_{i_1=0}^{j_{\sss 1}-1}p^{\walkvec{0}{1}\circ \walkvec{1}{i_1}}\left(\walkcoor{1}{i_1},\walkcoor{1}{i_1+1}\right)
\dots
\sum_{j_{\sss N}=0}^{\infty}\sum_{\walkvec{N}{j_{\sss N}+1}}|\Delta_{\sN}|\prod_{i_N=0}^{j_{\sss N}-1}p^{\walkvec{N-1}{j_{\sss {N-1}}+1}\circ \walkvec{N}{i_N}}\left(\walkcoor{N}{i_N},\walkcoor{N}{i_N+1}\right),\lbeq{absvalpi}
\enalign
where the sums over $j_{\sss k}$, $k\ge 2$ are all from $0$ to $\infty$.  Note that by \refeq{deltanbd}, $\Delta_{\sss 1}\neq 0$ only when $\eta^{(1)}_{j_1}=\eta^{(0)}_0=o$, and in particular, only when $j_1$ is odd, which is why $j_1$ is summed from $1$ onward.
We will use Lemma \ref{lem:togreens} to successively bound the sums over $j_{\sss k}$ in \refeq{absvalpi}, beginning with $k=N$.

When $N=1$, \refeq{absvalpi} becomes
\eqalign
&\sum_{\walkcoor{0}{1}}p^{o}(o,\walkcoor{0}{1})
\sum_{j_{\sss 1}=1}^{\infty}\sum_{\walkvec{1}{j_{\sss 1}+1}}
|\Delta_{\sss 1}|\prod_{i_1=0}^{j_{\sss 1}-1}p^{\walkvec{0}{1}\circ \walkvec{1}{i_1}}\left(\walkcoor{1}{i_1},\walkcoor{1}{i_1+1}\right) \nn \\
\le &\sum_{\walkcoor{0}{1}}p^{o}(o,\walkcoor{0}{1})
\sum_{j_{\sss 1}=1}^{\infty}\sum_{\walkvec{1}{j_{\sss 1}}}
P^{\walkvec{0}{1}}(\vec X_{j_{\sss 1}}=\vec\eta_{j_{\sss 1}})1_{\{\walkcoor{1}{j_1}=o\}}\epsilon_{\delta} \nn \\
=&\epsilon_{\delta}\sum_{\walkcoor{0}{1}}p^{o}(o,\walkcoor{0}{1})
\sum_{j_{\sss 1}=1}^{\infty}P^{\walkvec{0}{1}}(X_{j_{\sss 1}}=o)=\epsilon_{\delta}\sum_{j=2}^{\infty}P^{o}(X_{j}=o)\leq \eps_\delta\delta^{-1}G_q,
\enalign
where we used
\refeq{deltanbd} in the first inequality, and the last inequality follows by setting $i=0$, $u=o$ and $\vec\eta_m=\{o\}$ in \refeq{Gbound1}.

For $N\ge 2$, as above we write
\eqalign
\sum_{j_{\sss N}=0}^{\infty}\sum_{\walkvec{N}{j_{\sss N}+1}}|\Delta_{\sN}|\prod_{i_N=0}^{j_{\sss N}-1}p^{\walkvec{N-1}{j_{\sss {N-1}}+1}\circ \walkvec{N}{i_N}}\left(\walkcoor{N}{i_N},\walkcoor{N}{i_N+1}\right)\le&\sum_{j_{\sss N}=0}^{\infty}\sum_{\walkvec{N}{j_{\sss N}}}P^{\walkvec{N-1}{j_{\sss {N-1}}+1}}(\vec X_{j_N}=\walkvec{N}{j_{\sss N}})1_{\{\walkcoor{N}{j_{\sss N}}\in \walkvec{N-1}{j_{N-1}}\}}\epsilon_{\delta}\nn\\
\le &\epsilon_{\delta}\sum_{i=0}^{j_{N-1}}\sum_{j_{\sss N}=0}^{\infty}P^{\walkvec{N-1}{j_{\sss {N-1}}+1}}(X_{j_N}=\walkcoor{N-1}{i})\nn\\
\le &(j_{N-1}+1)\epsilon_{\delta}\delta^{-1}G_q.\lbeq{Npiece}
\enalign
For the sum over $j_{N-1}$, we proceed as above except that we now have an extra factor of $(j_{N-1}+1)$, whence we use \refeq{Gbound1} with $i=1$.  Continuing in this way until reaching the sum over $j_1$, we get $N-2$ factors of $\alpha=\epsilon_{\delta}\delta^{-2}G^{*2}_q$. For the sum over $j_1$, proceeding as for the $N=1$ case but with the extra factor $(j_1+1)$, we then have to deal with the term
\eqalign
&\sum_{\walkcoor{0}{1}}p^{o}(o,\walkcoor{0}{1})
\sum_{j_{\sss 1}=1}^{\infty}(j_1+1)\sum_{\walkvec{1}{j_{\sss 1}+1}}
|\Delta_{\sss 1}|\prod_{i_1=0}^{j_{\sss 1}-1}p^{\walkvec{0}{1}\circ \walkvec{1}{i_1}}\left(\walkcoor{1}{i_1},\walkcoor{1}{i_1+1}\right) \nn \\
\le &\epsilon_{\delta}\sum_{\walkcoor{0}{1}}p^{o}(o,\walkcoor{0}{1})
\sum_{j_{\sss 1}=1}^{\infty}(j_1+1)P^{\walkvec{0}{1}}(X_{j_{\sss 1}}=o)=\epsilon_{\delta}\sum_{j=2}^{\infty}jP^{o}(X_{j}=o)\le \epsilon_{\delta}\delta^{-2}G^{*2}_{q}=\alpha,
\enalign
where we have again applied \refeq{Gbound1}. Combining all the factors then gives \refeq{pibound}.
\qed
\vskip0.5cm

\subsection{The derivative of the speed formula}
\label{sec-derivpi}
From \refeq{speedformula2} and \refeq{mainderiv} we have that
\eq
\lbeq{deriv_form}
\frac{\del v^{[1]}}{\del \beta}= \kappa\rho+\frac{\del }{\del \beta}\sum_{m=2}^{\infty}\sum_{N=1}^{m-1}\sum_{x,y}(y^{[1]}-x^{[1]})\pi_m^{\sss(N)}(x,y),
\en
assuming that the latter derivative actually exists.

Recall \refeq{piNxydef} and define
\eqalign
\lbeq{phidef}
\varphi_m^{\smallsup{N}}(x,y):=&\frac{\del}{\del \beta} \pi_m^{\smallsup{N}}(x,y),
\enalign
which is well-defined as a finite sum of finite products of transition probabilities (see \refeq{trans_deriv}), and is non-zero only for a finite
set of $x,y\in\Z^d$ due to the nearest-neighbour constraint. In order to prove Theorem \ref{thm:monotone}, it is sufficient to show that
\eq
\lbeq{Thm1.3cond1}
\sup_{\beta\in [0,1]} \sum_{m=2}^{\infty}\sum_{N=1}^{m-1} |\sum_{x,y}  (y-x)^{[1]}\varphi_m^{\smallsup{N}}(x,y)|<\kappa\rho,
\en
and
\eq
\lbeq{Thm1.3cond2}
\lim_{m_0\uparrow\infty} \sup_{\beta\in [0,1]} \sum_{m=m_0}^\infty \sum_{N=1}^{m-1} |\sum_{x,y}  (y-x)^{[1]}\varphi_m^{\smallsup{N}}(x,y)|=0.
\en
This implies the uniform convergence (in $\beta\in [0,1]$) of the derivative of the summands in \refeq{deriv_form} and allows us to pass the derivative
inside the infinite sum and conclude that $\frac{\del v^{[1]}}{\del \beta}>0$.
\bigskip

Let us write
    \eqalign \varphi_m^{\smallsup{N}}(x,y)=\varphi_m^{\smallsup{N,1}}(x,y)+\varphi_m^{\smallsup{N,2}}(x,y)+\varphi_m^{\smallsup{N,3}}(x,y),\lbeq{varphibreak}
    \enalign
where by the product rule, $\varphi_m^{\smallsup{N,1}}(x,y)$, $\varphi_m^{\smallsup{N,2}}(x,y)$ and $\varphi_m^{\smallsup{N,3}}(x,y)$ arise from differentiating
$p^{o}(o,\walkcoor{0}{1})$, $\prod_{n=1}^{N}\prod_{i_{n}=0}^{j_{n}-1}p^{\walkvec{n-1}{j_{n-1}+1}\circ \walkvec{n}{i_{n}}}\left(\walkcoor{n}{i_{n}},\walkcoor{n}{i_{n}+1}\right)$ and $\prod_{n=1}^{N}\Delta_{\sss n}$, respectively, in \refeq{piNxydef}.

%Recall that in the current setting
%    \eqalign
%    \frac{\del}{\del\beta}{p_{\beta}^{\vec{\eta}_m}}(\eta_m,\eta_m+u)=
    %d^{-1}\big[1_{\{u=e_1\}}-1_{\{u=-e_1\}}\big]1_{\big\{\ell^{\vec{\eta}_m}(\eta_m,\eta_m+e_1)+\ell^{\vec{\eta}_m}(\eta_m,\eta_m-e_1)=0\big\}}
%   \enalign
%   which is bounded in absolute value by
%   \eqalign
%   d^{-1}1_{\{u=\pm e_1\}}\lbeq{pderivbound1}.
%   \enalign
%Similarly,
%\eqalign
%\frac{\del}{\del \beta}\Delta_{\sss i}(\walkcoor{i}{j_i},\walkcoor{i}{j_i}+u)=&d^{-1}\big[1_{\{u=-e_1\}}-1_{\{u=e_1\}}\big]
%1_{\big\{\ell^{\vec{\omega}^{\sss (i-1)}_{j_{i-1}+1}\circ \vec{\omega}^{\sss (i)}_{j_i}}(\walkcoor{i}{j_i},\walkcoor{i}{j_i}+e_1)>0, \ell^{\vec{\omega}^{\sss (i-1)}_{j_{i-1}+1}\circ \vec{\omega}^{\sss (i)}_{j_i}}(\walkcoor{i}{j_i},\walkcoor{i}{j_i}-e_1)=0\big\}} \nn \\
%&\times 1_{\big\{\ell^{\vec{\omega}^{\sss (i)}_{j_i}}(\walkcoor{i}{j_i},\walkcoor{i}{j_i}+e_1)+\ell^{\vec{\omega}^{\sss (i)}_{j_i}}(\walkcoor{i}{j_i},\walkcoor{i}{j_i}-e_1)=0\big\}},
%\enalign
%is bounded in absolute value by
%\eqalign
%d^{-1}1_{\{u=\pm e_1\}}
%1_{\big\{(\walkcoor{i}{j_i},\walkcoor{i}{j_i+1})\in \{(\walkcoor{i-1}{l},\walkcoor{i-1}{l+1}):0\le l\le j_{i-1}\}\big\}}.\nn \\
%\enalign
Define
\eqalign
F^{\smallsup{N}}:=&\sum_{\walkcoor{0}{1}}\big|\frac{\del}{\del \beta} p^{o}(o,\walkcoor{0}{1})\big|
\sum_{j_{\sss 1}=1}^{\infty}\sum_{\walkvec{1}{j_{\sss 1}+1}}
|\Delta_{\sss 1}|\prod_{i_1=0}^{j_{\sss 1}-1}p^{\walkvec{0}{1}\circ \walkvec{1}{i_1}}\left(\walkcoor{1}{i_1},\walkcoor{1}{i_1+1}\right)
\dots\nn\\
&\dots\sum_{j_{\sss N}=0}^{\infty}\sum_{\walkvec{N}{j_{\sss N}}}\prod_{i_N=0}^{j_{\sss N}-1}p^{\walkvec{N-1}{j_{\sss {N-1}}+1}\circ \walkvec{N}{i_N}}\left(\walkcoor{N}{i_N},\walkcoor{N}{i_N+1}\right)\big|\sum_{\walkcoor{N}{j_{N}+1}}(\walkcoor{N }{j_{N}+1}-\walkcoor{N}{j_{N}})^{[1]}\Delta_{\sN}\big|.\lbeq{rhoN}
\enalign
It follows that
\eq
\sum_{m=2}^\infty|\sum_{x,y}  (y-x)^{[1]}\varphi_m^{\smallsup{N,1}}(x,y)|\le F^{\smallsup{N}}.
\en
Similarly, for $k=1, \ldots, N$, let $H^{\smallsup{N}}_k$
be defined by replacing (in the definition  \refeq{piNxydef}):
$\Delta_{\sss n}$ with  $|\Delta_{\sss n}|$ for all $n=1, \ldots, N-1$, $\sum_{\eta^{(N)}_{j_N+1}}\Delta_{\sN}$ with $|\sum_{\eta^{(N)}_{j_N+1}}(\walkcoor{N }{j_{N}+1}-\walkcoor{N}{j_{N}})^{[1]}\Delta_{\sN}|$, $\prod_{i_{k}=0}^{j_{k}-1}p^{\walkvec{k-1}{j_{k-1}+1}\circ \walkvec{k}{i_{k}}}\left(\walkcoor{k}{i_{k}},\walkcoor{k}{i_{k}+1}\right)$ with
\eq
\sum_{l=0}^{j_{k}-1}\left|\frac{\del}{\del\beta}p^{\walkvec{k-1}{j_{k-1}+1}\circ \walkvec{k}{l}}\left(\walkcoor{k}{l},\walkcoor{k}{l+1}\right)\right|\prodtwo{i_{k}=0}{i_{k}\ne l}^{j_{k}-1}p^{\walkvec{k-1}{j_{k-1}+1}\circ \walkvec{k}{i_{k}}}\left(\walkcoor{k}{i_{k}},\walkcoor{k}{i_{k}+1}\right).
\en
For $k=1,\ldots, N$, let $J^{\smallsup{N}}_k$ be defined by replacing in \refeq{piNxydef}: $\Delta_{\sss i}$ with $|\Delta_{\sss i}|$
for all $1\leq i\leq N-1$ and $i\neq k$, $\Delta_{\sss k}$ with $|\frac{\del}{\del \beta}\Delta_{\sss k}|$ and $\sum_{\eta^{(N)}_{j_N+1}}\Delta_{\sN}$ with $|\sum_{\eta^{(N)}_{j_N+1}}(\walkcoor{N }{j_{N}+1}-\walkcoor{N}{j_{N}})^{[1]}\Delta_{\sN}|$ for $1\leq k\leq N-1$, and $\sum_{\eta^{(N)}_{j_N+1}}\Delta_{\sN}$
with $\big|\frac{\del}{\del \beta}\sum_{\walkcoor{N}{j_{N}+1}}(\walkcoor{N }{j_{N}+1}-\walkcoor{N}{j_{N}})^{[1]}\Delta_{\sN}\big|$ for $k=N$.

Letting $H^{\smallsup{N}}:=\sum_{k=1}^{N}H^{\smallsup{N}}_k$ and
$J^{\smallsup{N}}:=\sum_{k=1}^{N}J^{\smallsup{N}}_k$, we observe that
    \eqalign
    \sum_{m=2}^\infty |\sum_{x,y}  (y-x)^{[1]}\varphi_m^{\smallsup{N,2}}(x,y)|\le
    H^{\smallsup{N}}, \quad \text{and  }\sum_{m=2}^\infty |\sum_{x,y}  (y-x)^{[1]}\varphi_m^{\smallsup{N,3}}(x,y)|\le
    J^{\smallsup{N}}.\lbeq{3terms}
    \enalign

The remainder of this section is devoted to establishing bounds on $F^{\smallsup{N}}$, $J^{\smallsup{N}}$, and $H^{\smallsup{N}}$.

\begin{LEM}{\bf (Bounds on $F^{\smallsup{N}}$)}
\label{lem:rhobound}
We have
\eqalign
\lbeq{rhobound}
F^{\smallsup{N}} \le & \begin{cases}
\kappa\rho\, \epsilon_\delta\,\delta^{-1}\,G_q, & N=1, \\
\kappa\rho\, \delta^{-1}\,G_q \alpha^{N-1}, & N\geq 2.
\end{cases}
\enalign
\end{LEM}
\noindent
{\bf Proof.} When $N=1$, we first use \refeq{trans_deriv} and \refeq{Delta_N_bound} to get

\eqalign
F^{\smallsup{1}}\le &\sum_{\walkcoor{0}{1}}\kappa|\nu_2(\walkcoor{0}{1})-\nu_1(\walkcoor{0}{1})|
\sum_{j_{\sss 1}=1}^{\infty}\sum_{\walkvec{1}{j_{\sss 1}}}
\prod_{i_1=1}^{j_{\sss 1}-1}p^{\walkvec{0}{1}\circ \walkvec{1}{i_1}}\left(\walkcoor{1}{i_1},\walkcoor{1}{i_1+1}\right)
1_{\{\walkcoor{1}{j_1}=o\}}\rho \nn \\
=&\kappa\rho\sum_{\walkcoor{0}{1}}|\nu_2(\walkcoor{0}{1})-\nu_1(\walkcoor{0}{1})|
\sum_{j_{\sss 1}=1}^{\infty}P^{\walkvec{0}{1}}(X_{j_1}=o)\le \kappa\rho2(1-\delta)\delta^{-1}G_q,\lbeq{rho_first_bound}
\enalign
where we have used the fact that $\nu_1,\nu_2\in \Mi_{\gamma}(\Z_{d_*})$ with $\gamma\leq 1-\delta$, and we applied \refeq{Gbound1}.
This gives the bound for $N=1$.
%, this is bounded by
%\eq
%\kappa\rho G_q(0)\sum_{\walkcoor{0}{1}}|\nu_2(\walkcoor{0}{1})-\nu_1(\walkcoor{0}{1})|\le 2\gamma\kappa\rho G_q(0),\lbeq{rho_first_bound}
%\en

For $N\ge 2$ we proceed similarly, first using \refeq{Delta_N_bound} in the form
\[\big|\sum_{\walkcoor{N}{j_{N}+1}}(\walkcoor{N }{j_{N}+1}-\walkcoor{N}{j_{N}})^{[1]}\Delta_{\sN}\big|\le \rho 1_{\{\walkcoor{N}{j_N}\in \walkvec{N-1}{j_{N-1}}\}},\]
and then proceeding as in the proof of Proposition \ref{prp:pibound}.  This involves using \refeq{Gbound1} with $i=0$ to deal with what remains inside the sum over $j_N$,
which gives overall a factor of $(j_{N-1}+1) \rho \delta^{-1} G_q$. We then repeatedly use \refeq{deltanbd} and \refeq{Gbound1} with $i=1$ for the remaining terms inside the sums over $j_{N-1}, j_{N-2}, \dots, j_2$ in that order.  As in Proposition \ref{prp:pibound}, this gives $N-2$ factors of $\alpha$.  Finally we are left to deal with a term of the form
\eq
\kappa \sum_{\walkcoor{0}{1}}|\nu_2(\walkcoor{0}{1})-\nu_1(\walkcoor{0}{1})|
\sum_{j_{\sss 1}=1}^{\infty}(j_{\sss 1}+1)P^{\walkvec{0}{1}}(X_{j_1}=o)\le \kappa 2(1-\delta) \delta^{-2}G_q^{*2} = \kappa\alpha,
\en
where we again used \refeq{Gbound1} with $i=1$.
\qed

\begin{LEM}{\bf (Bounds on $J^{\smallsup{N}}$)}
\label{lem:chibound}
We have
\eqalign
\lbeq{Jbound}
J^{\smallsup{N}}_k \le & \begin{cases}
\kappa\rho\,\delta^{-2}(G_q(o)-\delta), & N=k=1, \\
\kappa\rho\,\delta^{-1}\,G_q \alpha^{N-1}, & N\geq 2, 1\leq k\leq N.
\end{cases}
\enalign
\end{LEM}
\noindent
{\bf Proof.} The second bound in \refeq{Jbound} follows essentially the same proof as that for Lemma \ref{lem:rhobound}. For $N\geq 2$, when we sum over $j_N$ and $\eta^{(N)}_{j_N+1}$, we apply \refeq{Delta_N_bound} if $k<N$ and we apply \refeq{Delta_N_deriv_bound} if $k=N$. When we sum over $j_k$ and $\eta^{(k)}_{j_k+1}$ with $k<N$, we apply \refeq{Delta_n_deriv_bound}.

To bound $J^{(1)}_1$, note that by \refeq{Delta_N_deriv_bound} applied to $\vec x_m=\vec\eta^{(0)}_1$ and $\vec \eta_n=\vec\eta^{(1)}_{j_1}$, we have
\eqalign
J^{\smallsup{1}}_1\le & \sum_{\walkcoor{0}{1}}p^o(o,\walkcoor{0}{1})
\sum_{j_{\sss 1}=1}^{\infty}\sum_{\walkvec{1}{j_{\sss 1}}}
\prod_{i_1=0}^{j_{\sss 1}-1}p^{\walkvec{0}{1}\circ \walkvec{1}{i_1}}\left(\walkcoor{1}{i_1},\walkcoor{1}{i_1+1}\right)
1_{\{\walkcoor{1}{j_1}=o\}}\kappa\rho=\kappa\rho\sum_{j=2}^{\infty}P^o(X_j=o).
\enalign
Using the parity of $\Z^{d_1}$, and proceeding as in the proof of Lemma \ref{lem:togreens},
\eqalign
\sum_{j=2}^{\infty}P^o(X_j=o)\le &\sum_{l=0}^{\infty}\mc{P}_q(Y_l=o)\sum_{j=l\vee 2}^{\infty}P^o(\mc{N}_j=l | Y_l=o)\nn\\
=&\sum_{j=2}^{\infty}P^o(\mc{N}_j=0)+\sum_{l=2}^{\infty}\mc{P}_q(Y_l=o)\sum_{j=l}^{\infty}P^o(\mc{N}_j=l |Y_l=0).
%=&E^o\big[\sum_{j=2}^{\infty}1_{\{\mc{N}_j=0\}}\big]+ (G_q(o)-1)E^o\big[\sum_{j=2}^{\infty}1_{\{\mc{N}_j=n\}}\big]
%\delta^{-1}\sum_{n=2}^{\infty}\mc{P}_q(X_n=o)=E^o\big[\sum_{j=2}^{\infty}1_{\{\mc{N}_j=0\}}\big]+\delta^{-1}[G_q(o)-1].\nn
\enalign
As shown in the proof of Lemma~\ref{lem:togreens}, for each $j\in\N$, conditional on $\vec{X}_{j-1}$ and $\vec Y_l$, either
$\mc{N}_j=\mc{N}_{j-1}$, or $\mc{N}_j=\mc{N}_{j-1}+1$ with probability at least $\delta$. Therefore by the same comparison
with the Bernoulli random walk $\mc{M}$ as in the proof of Lemma \ref{lem:qpart2}, we have
$$
\sum_{j=l}^{\infty}P^o(\mc{N}_j=l | Y_l=o)=E^o[\tau_{l+1}-\tau_l]\le \delta^{-1}.
$$
%is the expected amount of time $\mc{N}$ spends at position $n$, which satisfies $$.
% we have
%\begin{eqnarray*}
%\sum_{j=n}^{\infty}P^o(\mc{N}_j=n) = {\ch E}^o\big[\sum_{j=0}^{\infty}1_{\{\mc{N}_j=n\}}\big] &=& \sum_{k=1}^\infty %P^o\big(\sum_{j=0}^{\infty}1_{\{\mc{N}_j=n\}}\geq k \big)\\
%&\leq& 1+ (1-\delta)+(1-\delta)^2+\cdots = \delta^{-1}.
%\end{eqnarray*}
Similarly
$$
\sum_{j=2}^{\infty}P^o(\mc{N}_j=0)\leq \delta^{-1}-1.
$$
%= \sum_{j=0}^{\infty}P^o(\mc{N}_j=0) - P^o(\mc{N}_0=0) - P^o(\mc{N}_1=0) \leq \delta^{-1}-1.
Combining all the above bounds gives
\eq
J^{\smallsup{1}}_1\le \kappa\rho \big(\delta^{-1}-1 +\delta^{-1}[G_q(o)-1]\big) = \kappa\rho\delta^{-2}(G_q(o)-\delta),
\en
as required.
\qed

\vspace{.5cm}

To bound $H^{\smallsup{N}}_{k}$, we need a new lemma of the form of Lemma \ref{lem:togreens} which accommodates the derivative of the transition probability for one
of the steps.
\begin{LEM}
\label{lem:togreens2} Assume the same conditions as in Lemma \ref{lem:togreens}. Recall the definition of $\nu_1, \nu_2$ from assumption {\bf (A2)},
and the definition of edge local time $\ell(\vec \eta, V_{d_*})$ from \refeq{edgeloctime}. Then for each $i \in \Z_+$ and $u\in \Z^d$, we have
\eqalign
&\sum_{j=1}^{\infty}\frac{(j+i)!}{j!}\sum_{s=0}^{j-1}\sum_{\vec{\eta}_{\sss s}}
P^{\vec{x}_{\sss m}}(\vec{X}_{s}=\vec{\eta}_{\sss s})1_{\{\ell(\vec{x}_{\sss m}\circ\vec{\eta}_{\sss s},V_{d_*})=0\}} \sum_{\eta_{\sss s+1}}\kappa|\nu_2(\eta_{\sss s+1}-\eta_{\sss s})-\nu_1(\eta_{\sss s+1}-\eta_{\sss s})|\nn\\
&\qquad \qquad  \times P^{\vec{x}_{\sss m}\circ\vec{\eta}_{\sss s+1}}(X_{j-s-1}=u)
\ \ \le  \ \ \eps_\delta\, \kappa(i+1)!\delta^{-(i+2)}G^{*(i+2)}_q. \lbeq{lem3.7eq0}
%\\
%&\sum_{j=2}^{\infty}\frac{(j+i)!}{j!}\sum_{l=1}^{j}\sum_{\vec{\eta}_{\sss l}}
%P^{\vec{x}_{\sss m}}(\vec{X}_{l}=\vec{\eta}_{\sss l})\sum_{\eta_{\sss l+1}}\kappa|\nu_2(\eta_{\sss l+1}-\eta_{\sss l})-\nu_1(\eta_{\sss l+1}-\eta_{\sss l})|1_{\{\ell(\vec{x}_{\sss m}\circ\vec{\eta}_{\sss l},V_{d_*})=0\}}\\ &\times\sum_{\vec{\eta}_{\sss [l+2,\dots,j]}}P^{\vec{x}_{\sss m}\circ\vec{\eta}_{\sss l+1}}(\vec{X}_{j-l-1}=\vec{\eta}_{\sss [l+2,\dots,j]}) 1_{\{\eta_j=u\}}\le 2\gamma \kappa(i+1)!\mc{E}_{i+1}
%    &\sum_{j=0}^{\infty}\frac{(j+i)!}{j!}P^{\vec{\eta}_m}(S_{j}=u)\le i!\delta^{-(i+1)}G^{*(i+1)}_q,\lbeq{Gbound1b}\\
%    &\sum_{j=1}^{\infty}\frac{(j+i)!}{j!}P^{\vec{\eta}_m}(S_{j}=u)
%    \le i!\mc{E}_i,\lbeq{Gbound2b}\\
%    \text{ for any }w\in V_{d_0}, \quad &\sum_{j=1}^{\infty}\sum_{l=0}^{j-1}\frac{(j+i)!}{j!}P^{\vec{\eta}_m}(S_{j}=u|S_{l+1}-S_l=w)
%    \le (i+1)!\delta^{-(i+2)}G_q^{*(i+2)}\lbeq{Gbound3b}.
    \enalign
\end{LEM}
\noindent
{\bf Proof.} Since $\nu_1$ and $\nu_2$ have disjoint supports $\mc{S}_1$, and $\mc{S}_2$, the left hand side of \refeq{lem3.7eq0} equals
\eqalign\lbeq{lem3.7eq1}
\sum_{r=1}^2\sum_{j=1}^{\infty}\frac{(j+i)!}{j!}\sum_{s=0}^{j-1}\sum_{\vec\eta_s}P^{\vec{x}_m}(\vec{X}_{s}=\vec\eta_s)1_{\{\ell(\vec{x}_{\sss m}\circ\vec{\eta}_{\sss s},V_{d_*})=0\}} \!\!\! \sum_{w:=\eta_{s+1}-\eta_s\in\mc{S}_r}\!\!\!\!\!\!\!\! \kappa\nu_r(w) P^{\vec{x}_{\sss m}\circ\vec{\eta}_{\sss s+1}}(X_{j-s-1}=u).
\enalign
Let $P_{(s,w)}^{\vec{x}_m}$ denote the law of a random walk which evolves according to a RWRE with history $\vec{x}_m$, except that the $s+1$-st step is deterministic
and equals $w$, and if this step has zero probability given the history of $X$ up to time $s$, then the walk is killed. More precisely,
\[
P_{(s,w)}^{\vec{x}_m}(X_{n+1}-X_n=z|\vec{X}_n=\vec{x}_n)=\begin{cases}
\delta_w(z), &\text{ if }n=s,\\
P^{\vec{x}_m}(X_{n+1}-X_n=z|\vec{X}_n=\vec{x}_n), &\text{ otherwise,}
\end{cases}\]
and the walk is killed at time $s+1$ if $P^{\vec{x}_m}(X_{s+1}-X_s=w|\vec{X}_s=\vec{x}_s)=0$.

Then \refeq{lem3.7eq1} is bounded by
\eqalign\lbeq{3.7E2}
\kappa\sum_{r=1}^2\sum_{w\in\mc{S}_r} \nu_r(w) \sum_{j=1}^{\infty}\frac{(j+i)!}{j!}\sum_{s=0}^{j-1}P_{(s,w)}^{\vec{x}_m}(X_j=u).
\enalign
As in the proof of Lemma \ref{lem:togreens}, we have
\eq
P_{(s,w)}^{\vec{x}_m}(X_j=u)\le \sum_{l=0}^{j-1}P_q(Y_l=\Pi_{d_1}(u))P_{(s,w)}^{\vec{x}_m}(\mc{N}_j=l|Y_l=\Pi_{d_1}(u)),
\en
where $\mc{N}_j$ is the number of steps of $\vec{X}_j$ with non-zero $\Pi_{d_1}$-projection.  Thus, the sum over $j$ in \refeq{3.7E2} is bounded by
\eqalign
&\sum_{j=1}^{\infty}\frac{(j+i)!}{j!}\sum_{s=0}^{j-1}\sum_{l=0}^{j-1}P_q(Y_l=\Pi_{d_1}(u))P_{(s,w)}^{\vec{x}_m}(\mc{N}_j=l|Y_l=\Pi_{d_1}(u))\nn\\
=&\sum_{l=0}^{\infty}P_q(Y_l=\Pi_{d_1}(u))\sum_{s=0}^{\infty}\sum_{j=(l\vee s)+1}^{\infty}\frac{(j+i)!}{j!}P_{(s,w)}^{\vec{x}_m}(\mc{N}_j=l|Y_l=\Pi_{d_1}(u)).\lbeq{banana1}
\enalign
As in the proof of Lemma \ref{lem:qpart2}, we can couple $\mc{N}$ with a Bernoulli random walk $\mc{M}$ such that
\eqalign
&\sum_{j=(l\vee s)+1}^{\infty}\frac{(j+i)!}{j!}P_{(s,w)}^{\vec{x}_m}(\mc{N}_j=l|Y_l=\Pi_{d_1}(u)) \nn \\
=&E_{(s,w)}^{\vec{x}_m}\Big[\sum_{j=\tau_l\vee (s+1)}^{\tau_{l+1}-1}\frac{(j+i)!}{j!}\big|Y_l=\Pi_{d_1}(u)\Big] \nn\\
\le&  E\Big[\sum_{j=\tau_l(\mc{M}-1)\vee (s+1)}^{\tau_{l+1}(\mc{M}-1)-1}\frac{(j+i)!}{j!}\Big]
=\sum_{j=(l\vee s)+1}^{\infty}\frac{(j+i)!}{j!}P(\mc{M}_j-1=l)
\enalign
where $\mc{M}_n\sim$Bin$(n,\delta)$.  Therefore the summation over $s$ in \refeq{banana1} is bounded by
\eqalign
& \sum_{s=0}^{\infty}\sum_{j=(l\vee s)+1}^{\infty}\frac{(j+i)!}{j!}P(\mc{M}_j-1=l) =\sum_{j=l+1}^{\infty}\sum_{s=0}^{j-1}\frac{(j+i)!}{j!}P(\mc{M}_j-1=l) \nn \\
=& \sum_{j=l+1}^{\infty}\frac{(j+i)!}{(j-1)!}P(\mc{M}_j=l+1) < \sum_{j=l}^{\infty}\frac{(j+i+1)!}{j!}P(\mc{M}_j=l) \leq \delta^{-(i+1)}\frac{(l+i+1)!}{l!},\lbeq{forend}
\enalign
as in \refeq{comb1}. Substituting this bound back into \refeq{banana1} and then into \refeq{3.7E2} then proves the lemma, where we need to use \refeq{Gform}
and the fact that $\nu_1, \nu_2$ each has total mass $\gamma \leq 1-\delta = \eps_\delta/2$.
\qed
\medskip

\begin{LEM}{\bf (Bounds on $H^{\smallsup{N}}$)}
\label{lem:gammabound}
We have
\[H^{\smallsup{N}}_{\sss k}\le \begin{cases}
\kappa\rho\,\alpha^{N}, & N=k\geq 1,\\
2\kappa\rho\, \eps_\delta^2 \delta^{-4} G_q G_q^{*3} \alpha^{N-2}, &N>k\geq 1.
\end{cases}\]
\end{LEM}
\noindent
{\bf Proof.} We first use \refeq{trans_deriv} and \refeq{Delta_N_bound} to get
\eqalign
H^{\smallsup{1}}_{\sss 1}\ \le\ &\ \sum_{\walkcoor{0}{1}}p^o(o,\walkcoor{0}{1})
\sum_{j_{\sss 1}=1}^{\infty}\sum_{l=0}^{j_{\sss 1}-1} \sum_{\walkvec{1}{l}}P^{\walkvec{0}{1}}(\vec{X}_{l}=\walkvec{1}{l})
1_{\{\ell(\walkvec{0}{1}\circ\walkvec{1}{l},V_{d_*})=0\}}
\nn \\
&\qquad \times \sum_{\walkcoor{1}{l+1}}\kappa|\nu_2(\walkcoor{1}{l+1}-\walkcoor{1}{l})-\nu_1(\walkcoor{1}{l+1}-\walkcoor{1}{l})|\,
P^{\walkvec{0}{1}\circ\walkvec{1}{l+1}}(X_{j_1-l-1}=o) \rho \nn\\
\ =\ &\ \sum_{j=2}^{\infty}\sum_{l=1}^{j}\sum_{\vec{\eta}_{\sss l}} p^{o}(\vec{X}_{l}=\vec{\eta}_{\sss l})\sum_{\eta_{\sss l+1}}\kappa|\nu_2(\eta_{\sss l+1}-\eta_{\sss l})-\nu_1(\eta_{\sss l+1}-\eta_{\sss l}) |1_{\{\ell(\vec{\eta}_{\sss l},V_{d_*})=0\}} P^{\vec{\eta}_{\sss l+1}}(X_{j-l-1}=0)\rho \nn\\
\ \le\ &\ \kappa \rho \eps_\delta \delta^{-2} G_q^{*2} = \kappa\rho\alpha,
\enalign
where we have used Lemma \ref{lem:togreens2} with $i=0$.

For $N\ge 2$, we first bound $|\sum_{\walkcoor{N}{j_{N}+1}}(\walkcoor{N }{j_{N}+1}-\walkcoor{N}{j_{N}})^{[1]}\Delta_{\sN}|$ in
$H^{\smallsup{N}}_{\sss k}$ by $\rho 1_{\{\eta^{(N)}_{j_N}\in \vec\eta^{(N-1)}_{j_{N-1}}\}}$ using \refeq{Delta_N_bound}.

For $k=N$, we then use Lemma \ref{lem:togreens2} with $i=0$ to bound the sum over $j_{\sss N}$, yielding a factor of $\kappa\rho\eps_\delta\delta^{-2}G_q^{*2}(j_{N-1}+1)$. We then proceed exactly as in the proof of Proposition \ref{prp:pibound} on the sums over $j_{\sss N-1},\dots, j_{\sss 1}$ in that order, giving $N-1$ factors of $\alpha$.

If $N>k>1$, we use \refeq{Gbound1} with $i=0$ to bound the sum over $j_{\sss N}$, yielding a factor of $\rho\delta^{-1}G_q(j_{N-1}+1)$. We then proceed exactly as in the proof of Proposition \ref{prp:pibound} on the sums over $j_{\sss N-1},\dots, j_{\sss k+1}$ in that order, giving $N-(k+1)$ factors of $\alpha$.  We then obtain a factor $\epsilon_{\delta}$ from $\sum_{\eta^{(k)}_{j_k+1}}|\Delta_{\sss k}|$ and then we use Lemma \ref{lem:togreens2} with $i=1$ to bound the sum over $j_k$, giving a factor $2\kappa\eps_\delta\delta^{-3}G^{*3}_q(j_{k-1}+1)$.  We then proceed exactly as in the proof of Proposition \ref{prp:pibound} on the sums over $j_{\sss k-1},\dots, j_{\sss 1}$ in that order giving $k-1$ additional factors of $\alpha$. The case $N>k=1$ is similar.
\qed

\section{Proof of Theorem \ref{thm:monotone}}
\label{sec:mono_proof}
By assumption {\bf (A3)}, $G^{*i}_q<\infty$ for $i=1,2,3,4$.  It follows that $\alpha=2(1-\delta)\delta^{-2}G^{*2}_q<1$ for $\delta<1$ sufficiently close to $1$ (depending only on $q(\cdot)$), in which case \refeq{pibound} is summable in $N$, and thus \refeq{speedformula} and \refeq{speedformula2} hold. By assumption {\bf (A3)} that $G_q(o)<2$, and by Lemmas \ref{lem:rhobound}, \ref{lem:chibound}, and \ref{lem:gammabound}, uniformly in $\beta\in [0,1]$, we have
\eqalign
&\sum_{N\ge 1}F^{\sss (N)}\le \kappa\rho C_{1,\delta},\quad \sum_{N\ge 1}\sum_{k=1}^N J^{\sss (N)}_k\le \kappa\rho C_{2,\delta},\quad \sum_{N\ge 1}\sum_{k=1}^NH^{\sss (N)}_k\le \kappa\rho C_{3,\delta},\lbeq{threebounds}
\enalign
where $C_{1,\delta}, C_{3,\delta}\searrow 0$ as $\delta\nearrow 1$ and $C_{2,\delta}<1$ for $\delta$ sufficiently close to 1. It follows that for $\delta$ sufficiently close to 1 (depending only on $q(\cdot)$), \refeq{Thm1.3cond1} holds. To conclude the proof of Theorem~\ref{thm:monotone}, it only remains to verify \refeq{Thm1.3cond2}.

Note that
$$
\sum_{m=1}^\infty \sum_{N=N_0}^{m-1} |\sum_{x,y}  (y-x)^{[1]}\varphi_m^{\smallsup{N}}(x,y)| \leq \sum_{N=N_0}^\infty( F^{\sss (N)}+J^{\sss (N)}+H^{\sss (N)}),
$$
which by Lemmas \ref{lem:rhobound}, \ref{lem:chibound}, and \ref{lem:gammabound} tends to $0$ (uniformly in $\beta\in [0, 1]$) as $N_0\to\infty$ if $\delta$ is
sufficiently close to $1$. Therefore \refeq{Thm1.3cond2} follows from Lemma \ref{lem:mtailbd} below, which concludes the proof of Theorem~\ref{thm:monotone}.
\qed

\begin{LEM}\label{lem:mtailbd}
If $\delta\in (0,1)$ is chosen such that $\alpha=2(1-\delta)\delta^{-2}G^{*2}_q<1$, then for each $N\in\N$,
\eq\lbeq{summbd}
\lim_{m_0\uparrow \infty}\sup_{\beta\in [0,1]}\sum_{m=m_0}^\infty |\sum_{x,y}  (y-x)^{[1]}\varphi_m^{\smallsup{N}}(x,y)|=0.
\en
\end{LEM}
\noindent
{\bf Proof.} As in \refeq{varphibreak}, we will split $\varphi_m^{\smallsup{N}}$ into $\varphi_m^{\smallsup{N, i}}$ for $1\leq i\leq 3$. It will be sufficient to
verify \refeq{summbd} with $\varphi_m^{\smallsup{N}}$ replaced by $\varphi_m^{\smallsup{N, i}}$ for $i=1,2,3$. Note that
\eqalign
\sum_{m=m_0}^\infty |\sum_{x,y}  (y-x)^{[1]}\varphi_m^{\smallsup{N,1}}(x,y)|
\leq & \sum_{j_1,\cdots, j_N\geq 0\atop j_1+\cdots+j_N \geq {m_0-N}}
\sum_{\walkcoor{0}{1}}\big|\frac{\del}{\del \beta} p^{o}(o,\walkcoor{0}{1})\big|
\sum_{\walkvec{1}{j_{\sss 1}+1}}  P^{\vec \eta^{(0)}_1}(\vec X_{j_1}=\vec\eta^{(1)}_{j_1})|\Delta_{\sss 1}|
\dots\nn\\
&\dots\sum_{\walkvec{N}{j_{\sss N}}} P^{\vec \eta^{(N-1)}_{j_{N-1}+1}}(\vec X_{j_N}=\vec\eta^{(N)}_{j_N})
\big|\sum_{\walkcoor{N}{j_{N}+1}}(\walkcoor{N }{j_{N}+1}-\walkcoor{N}{j_{N}})^{[1]}\Delta_{\sN}\big| \nn \\
\leq & \sum_{k=1}^N \sum_{{i\neq k: j_i\geq 0 \atop j_k\geq \frac{m_0}{N}-1}} \sum_{\walkcoor{0}{1}}\big|\frac{\del}{\del \beta} p^{o}(o,\walkcoor{0}{1})\big|
\sum_{\walkvec{1}{j_{\sss 1}+1}}  P^{\vec \eta^{(0)}_1}(\vec X_{j_1}=\vec\eta^{(1)}_{j_1})|\Delta_{\sss 1}|
\dots\nn\\
&\dots\sum_{\walkvec{N}{j_{\sss N}}} P^{\vec \eta^{(N-1)}_{j_{N-1}+1}}(\vec X_{j_N}=\vec\eta^{(N)}_{j_N})
\big|\sum_{\walkcoor{N}{j_{N}+1}}(\walkcoor{N }{j_{N}+1}-\walkcoor{N}{j_{N}})^{[1]}\Delta_{\sN}\big|, \lbeq{sumoverk}
\enalign
where we made the observation that one of the $N$ paths $\vec\eta_{j_N}$ must have length at least $\frac{m_0}{N}-1$. In the sum over $k$ in \refeq{sumoverk}, if $k=N$,
then the sum over $j_N\geq \frac{m_0}{N}-1$ can be bounded by
\eq\lbeq{tb1}
\rho \sum_{{j_N\geq \frac{m_0}{N}-1}} \sum_{i=0}^{j_{N-1}} P^{\vec \eta^{(N-1)}_{j_{N-1}+1}}(\vec X_{j_N}=\vec\eta^{(N-1)}_{i})
\leq \rho (j_{N-1}+1) \sup_{u\in\Z^d} \sum_{{j_N\geq \frac{m_0}{N}-1}}  P^{\vec \eta^{(N-1)}_{j_{N-1}+1}}(\vec X_{j_N}=u).
\en
By Lemma \ref{lem:togreens}, for each $i\geq 0$, we have
\eq\lbeq{chebybd1}
\sup_{u\in\Z^d}  \sum_{j=m_0}^{\infty}\frac{(j+i)!}{j!}P^{\vec{\eta}_m}(X_{j}=u) \le  \sup_{u\in\Z^d} \sum_{j=0}^{\infty}\frac{(j+i+1)!}{m_0\, j! }P^{\vec{\eta}_m}(X_{j}=u) \leq \frac{(i+1)!\delta^{-(i+2)}G^{*(i+2)}_q}{m_0}.
\en

Applying this bound with $i=0$ to \refeq{tb1} then gives a factor of $1/m_0$. Summing over $j_{N-1},\cdots, j_1$ in \refeq{sumoverk} with $k=N$
as in the proof of Lemma \ref{lem:rhobound} then gives a bound proportional to $1/m_0$, which is independent of $\beta$ and tends to $0$ as $m_0\to\infty$.

In the sum over $k$ in \refeq{sumoverk}, if $k<N$, then we sum over $j_N, \cdots, j_{k+1}$ similarly as in the proof of Lemma \ref{lem:rhobound}, each sum giving rise
to a constant factor depending only on $\delta$ and $q(\cdot)$. When we sum over $j_k\geq \frac{m_0}{N}-1$, we need to bound
$$
\sum_{{j_k\geq \frac{m_0}{N}-1}} (j_k+1)\sum_{i=0}^{j_{k-1}} P^{\vec \eta^{(k-1)}_{j_{k-1}+1}}(\vec X_{j_k}=\vec\eta^{(k-1)}_{i}),
$$
for which we can apply \refeq{chebybd1} with $i=1$ to obtain a factor of $1/m_0$. Summing over $j_{k-1},\cdots, j_1$ only leads to bounded constant factors.
This verifies \refeq{summbd} with $\varphi_m^{\smallsup{N,1}}$ in place of $\varphi_m^{\smallsup{N}}$.

The proof of \refeq{summbd} with $\varphi_m^{\smallsup{N}}$ replaced by $\varphi_m^{\smallsup{N,2}}$ or $\varphi_m^{\smallsup{N,3}}$ is similar. Of the $N$ paths $\vec\eta^{(k)}_{j_k+1}$, $1\leq k\leq N$, one of these will have length at least ${j_k\geq \frac{m_0}{N}-1}$. We then draw upon the proofs of Lemmas \ref{lem:chibound} and \ref{lem:gammabound} to sum over $j_k$, $1\leq k\leq N$, and apart from \refeq{chebybd1}, we will also need the following corollary of Lemma \ref{lem:togreens2}:
\eqalign
&\sup_{u\in\Z^d} \sum_{j=m_0}^{\infty}\frac{(j+i)!}{j!}\sum_{l=0}^{j-1}\sum_{\vec{\eta}_{\sss l}}
P^{\vec{x}_{\sss m}}(\vec{X}_{l}=\vec{\eta}_{\sss l})1_{\{\ell(\vec{x}_{\sss m}\circ\vec{\eta}_{\sss l},V_{d_*})=0\}} \sum_{\eta_{\sss l+1}}\kappa|\nu_2(\eta_{\sss l+1}-\eta_{\sss l})-\nu_1(\eta_{\sss l+1}-\eta_{\sss l})|\nn\\
&\qquad \qquad  \times P^{\vec{x}_{\sss m}\circ\vec{\eta}_{\sss l+1}}(X_{j-l-1}=u)
\ \ \le  \ \ \frac{\eps_\delta\, \kappa(i+2)!\delta^{-(i+3)}G^{*(i+3)}_q}{m_0}. \lbeq{4.7}
    \enalign
Note that we would need to apply \refeq{4.7} with $i=1$, which uses $G^{*4}_q<\infty$ from assumption {\bf (A3)}. The details will be left to the reader.
\qed
\bigskip

\subsection{Remark on Assumption (A3) in Theorem \ref{thm:monotone}}
%\subsection{Remark on (b) versus (b') in Theorem \ref{thm:monotone}}
\label{sec:G4}
It is possible to replace the assumption that $G^{*4}<\infty$ with a local central limit theorem bound of the form
\eq
 \sup_x \mc{P}_q(X_n=x)\le \frac{C_q}{n^a}, \quad \text{for some }a>3.\lbeq{localCLT}
\en
The usual choice would be $a=\frac{d_1}{2}$, which is greater than 3 when $d_1>6$.   We give the main ideas of the argument here.  For further details see \cite{H08compete}.

The bound $G^{*4}<\infty$ was used when obtaining the estimate \refeq{4.7}.  However we only require that the left hand side of \refeq{4.7} converges to 0 as $m_0\ra\infty$.  As in the proof of Lemma \ref{lem:togreens2} (c.f.~\refeq{banana1}) this involves estimating
\eqalign\lbeq{end1}
&\sup_u\kappa\sum_{r=1}^2\sum_{w\in\mc{S}_r} \nu_r(w) \sum_{j=m_0}^{\infty}\frac{(j+i)!}{j!}\sum_{s=0}^{j-1}P_{(s,w)}^{\vec{x}_m}(X_j=u)\nn\\
\le &\sup_u \kappa\sum_{r=1}^2\sum_{w\in\mc{S}_r} \nu_r(w)\sum_{l=0}^{\infty}P_q(Y_l=\Pi_{d_1}(u))\sum_{s=0}^{\infty}\sum_{j=l\vee s\vee m_0}^{\infty}\frac{(j+i)!}{j!}P_{(s,w)}^{\vec{x}_m}(\mc{N}_j=l|Y_l=\Pi_{d_1}(u))\nn\\
\le &\sup_u\kappa\sum_{r=1}^2\sum_{w\in\mc{S}_r} \nu_r(w)\sum_{l=0}^{\infty}P_q(Y_l=\Pi_{d_1}(u))\sum_{s=0}^{\infty}\sum_{j=l\vee s\vee m_0}^{\infty}\frac{(j+i)!}{j!}P(\mc{M}_j-1=l).
\enalign
Using the local CLT bound \refeq{localCLT} this is bounded by
\eqalign
\eps_{\delta}\kappa \sum_{l=0}^{\infty}\frac{C_q}{l^a}\sum_{j=l\vee m_0}^{\infty}\frac{(j+i+1)!}{j!}P(\mc{M}_j-1=l)\le & \eps_{\delta}\kappa\sum_{l=0}^{K}\frac{C_q}{l^a}\sum_{j=l\vee m_0}^{\infty}\frac{(j+i+1)!}{j!}P(\mc{M}_j-1=l)\nn\\
&+\eps_{\delta}\kappa\sum_{l=K}^{\infty}\frac{C_q}{l^a}\sum_{j=l}^{\infty}\frac{(j+i+1)!}{j!}P(\mc{M}_j-1=l).\lbeq{end2}
\enalign
As in \refeq{forend} the last term of \refeq{end2} is bounded by
\eq
\eps_{\delta}\kappa\sum_{l=K}^{\infty}\frac{C_q}{l^a}\delta^{-(i+1)}\frac{(l+i+1)!}{l!}\le \eps_{\delta}\kappa\delta^{-(i+1)}\sum_{l=K}^{\infty}\frac{C'_q}{l^{a-(i+1)}}\nn,
\en
which can be made arbitrarily small by choosing $K$ large depending on $q$ and $i\in \{0,1\}$ when $a-(i+1)>1$, (i.e.~$a>3$ when $i=1$).
For any $K$, the first term of \refeq{end2} can be made arbitrarily small by choosing $m_0$ sufficiently large (see the proof of Lemma 4.1 of \cite{H08compete} for details).

\subsection*{Acknowledgements}
MH thanks Remco van der Hofstad and Alain-Sol Sznitman for helpful discussions at the initial stages of this project, and NUS mathematics department for support during his visit in 2009. We thank Jonathon Peterson for pointing out a false statement in an earlier version, and an anonymous referee for pointing out an error in our previous formulation of Theorem \ref{T:LLN}.

\end{document}